\documentclass[preprint,12pt]{elsarticle}

\usepackage{amssymb}
\newtheorem{thm}{Theorem}[section]

\newdefinition{defin}{Definition}[section]
\newproof{pf}{{\bf Proof}}

\makeatletter
\def\subsection{\@startsection{subsection}{2}{\z@}{.0ex plus-0ex
    minus-.0ex}{.0ex plus.0ex}{\reset@font\indent\bf}}
\makeatother

\newcommand{\vep}{\varepsilon}
\newcommand{\la}{\lambda}
\newcommand{\vfi}{\varphi}
\newcommand{\NB}{\mathbb{N}}  \newcommand{\RB}{\mathbb{R}}
 \newcommand{\XB}{\mathbb{X}}
\newcommand{\D}{\displaystyle}
\newcommand{\ov}[1]{\overline{#1}}
\newcommand{\eq}[1]{{\rm(\ref{#1})}}
\newcommand{\GV}[1]{\mathop{\rm GV}_{\!#1}\nolimits}
\newcommand{\wto}{\stackrel{\!\!w}{\to}}
\newcommand{\is}{i${}_{\mbox{\footnotesize\rm s}}$}

\catcode`\@=11  \@addtoreset{equation}{section}  \catcode`\@=12

\journal{arXiv}

\begin{document}

\begin{frontmatter}

\title{A fixed point theorem for contractions\\
   in modular metric spaces\tnoteref{t1}}
\tnotetext[t1]{Supported by Scientific Foundation of the National Research University
Higher School of Economics, Individual Research Project No.~10-01-0071.}

\author[hsenn]{Vyacheslav V.~Chistyakov}%\corref{cor1}}
\ead{czeslaw@mail.ru, vchistyakov@hse.ru}
\address[hsenn]{Department of Applied Mathematics and Computer Science,
National Research University Higher School of Economics, Bol'shaya Pech{\"e}rskaya
Street 25/12,\\ Nizhny Novgorod 603155, Russian Federation}

\begin{abstract}
The notion of a (metric) {\em modular\/} on an arbitrary set and the corresponding
{\em modular space}, more general than a metric space, were introduced and studied
recently by the author [V.\,V.~Chistyakov, Metric modulars and their application,
Dokl.\ Math.\ 73\,(1) (2006) 32--35, and Modular metric spaces, I: Basic concepts,
Nonlinear Anal.\ 72\,(1) (2010) 1--14].
In this paper we establish a fixed point theorem for contractive maps in modular spaces.
It is related to contracting rather ``generalized average velocities'' than metric
distances, and the successive approximations of fixed points converge to the fixed points
in a weaker sense as compared to the metric convergence.
\end{abstract}


\end{frontmatter}

\section{Introduction} \label{s:intro}

The metric fixed point theory (\cite{GoebelKirk,KirkSims}) and its variations
(\cite{HadzicPap}) are far-reaching developments of Banach's Contraction Principle,
where {\em metric conditions\/} on the underlying space and maps under
consideration play a fundamental role. This paper addresses fixed points of
nonlinear maps in {\em modular spaces\/} introduced recently by the author
(\cite{posit}--\cite{NA-II}) as generalizations of Orlicz spaces and classical modular
spaces (\cite{KR,Maligranda}, \cite{Musielak}--\cite{RaoRen}), where {\em modular
structures\/} (involving nonlinearities with more rapid growth than power-like functions),
play the crucial role. Under different contractive assumptions and the supplementary
$\Delta_2$-condition on modulars fixed point theorems in classical modular linear
spaces were established in~\cite{ait,Kh-Ko-Re,Kh-Ko-Sh}.

\smallbreak
We begin with a certain motivation of the definition of a (metric) {\em modular},
introduced axiomatically in \cite{DAN06,NA-I}. A simple and natural way to do it
is to turn to physical interpretations. Informally speaking, whereas a metric on a set
represents nonnegative finite distances between any two points of the set,
a modular on a set attributes a nonnegative (possibly, infinite valued) ``field of (generalized)
velocities'': to each ``time'' $\lambda\!>\!0$ (the absolute value of) an average velocity
$w_\lambda(x,y)$ is associated in such a way that in order to cover the ``distance''
between points $x,y\in X$ it takes time $\lambda$ to move
from $x$ to $y$ with velocity $w_\lambda(x,y)$. Let us comment on this
in more detail by exhibiting an appropriate example. If $d(x,y)\!\ge\!0$ is the distance
from $x$ to $y$ and a number $\lambda\!>\!0$ is interpreted as time, then the~value 
  \begin{equation} \label{e:dxy}
w_\lambda(x,y)=\frac{d(x,y)}\lambda
  \end{equation}
is the average velocity, with which one should move from $x$ to $y$ during time
$\lambda$, in order to cover the distance $d(x,y)$.
The following properties of the quantity from \eq{e:dxy} are quite natural.

(i) Two points $x$ and $y$ from $X$ coincide (and $d(x,y)=0$) if and only if
any time $\lambda>0$ will do to move from $x$ to $y$ with velocity $w_\lambda(x,y)=0$
(i.e., no movement is needed at any time). Formally, given $x,y\in X$, we have:
  \begin{equation} \label{e:axi}
\mbox{$x=y$ \,iff \,$w_\lambda(x,y)=0$ \,for \,all \,$\lambda>0$ \,(nondegeneracy),}
  \end{equation}
where `iff' means as usual `if and only if'.

(ii) Assuming the distance function to be symmetric, $d(x,y)=d(y,x)$, we find that
for any time $\lambda>0$ the average velocity during the movement from $x$ to $y$
is the same as the average velocity in the opposite direction, i.e.,
for any $x,y\in X$ we have:
  \begin{equation} \label{e:axi2}
\mbox{$w_\lambda(x,y)=w_\lambda(y,x)$ \,for \,all \,$\lambda>0$ \,(symmetry).}
  \end{equation}

(iii) The third property of \eq{e:dxy}, which is, in a sense, a counterpart of the
triangle inequality (for velocities!), is the most important. Suppose the movement
from $x$ to $y$ happens to be made in two different ways, but the {\em duration
of time is the same\/} in each case: (a) passing through a third point $z\in X$, or
(b) straightforward from $x$ to $y$. If $\lambda$ is the time needed to get from
$x$ to $z$ and $\mu$ is the time needed to get from $z$ to $y$, then the corresponding
average velocities are $w_\lambda(x,z)$ (during the movement from $x$ to $z$) and
$w_\mu(z,y)$ (during the movement from $z$ to $y$). The total time needed for
the movement in the case (a) is equal to $\lambda+\mu$. Thus, in order to move
from $x$ to $y$ as in the case (b) one has to have the average velocity equal to
$w_{\lambda+\mu}(x,y)$. Since (as a rule) the straightforward distance $d(x,y)$
does not exceed the sum of the distances $d(x,z)+d(z,y)$, it becomes clear from
the physical intuition that the velocity $w_{\lambda+\mu}(x,y)$ does not exceed
at least one of the velocities $w_\lambda(x,z)$ or $w_\mu(z,y)$.
Formally, this is expressed as
  \begin{equation} \label{e:axi3}
w_{\lambda+\mu}(x,y)\le\max\{w_\lambda(x,z),w_\mu(z,y)\}\le
w_\lambda(x,z)+w_\mu(z,y)
  \end{equation}
for all points $x,y,z\in X$ and all times $\lambda,\mu>0$ (``triangle'' inequality). In fact,
these inequalities can be verified rigorously: if, on the contrary, we assume that
$w_\lambda(x,z)<w_{\lambda+\mu}(x,y)$ and $w_\mu(z,y)<w_{\lambda+\mu}(x,y)$,
then multiplying the first inequality by $\lambda$, the second inequality---by $\mu$,
summing the results and taking into account \eq{e:dxy}, we find
  $d(x,z)=\lambda w_\lambda(x,z)<\lambda w_{\lambda+\mu}(x,y)$
and
  $d(z,y)=\mu w_{\mu}(z,y)<\mu w_{\lambda+\mu}(x,y),$
and it follows that
  $d(x,z)+d(z,y)<(\lambda+\mu)w_{\lambda+\mu}(x,y)=d(x,y),$
which contradicts the triangle inequality for~$d$.

Inequality \eq{e:axi3} can be obtained in a little bit more general situation.
Let $f:(0,\infty)\to(0,\infty)$ be a function from the set of positive reals into
itself such that the function $\lambda\mapsto\lambda/f(\lambda)$ is nonincreasing
on $(0,\infty)$. Setting $w_\lambda(x,y)=d(x,y)/f(\lambda)$ (note that
$f(\lambda)=\lambda$ in \eq{e:dxy}), we have
  \begin{eqnarray}
w_{\lambda+\mu}(x,y)\!\!&=\!\!&\frac{d(x,y)}{f(\lambda\!+\!\mu)}\le
  \frac{d(x,z)\!+\!d(z,y)}{f(\lambda\!+\!\mu)}\le
  \frac{\lambda}{\lambda\!+\!\mu}\!\cdot\!\frac{d(x,z)}{f(\lambda)}+
  \frac{\mu}{\lambda\!+\!\mu}\!\cdot\!\frac{d(z,y)}{f(\mu)}\le\nonumber\\
\!\!&\le\!\!&\frac{\lambda}{\lambda+\mu}\,w_\lambda(x,z)
  +\frac{\mu}{\lambda+\mu}\,w_\mu(z,y)\le
  w_\lambda(x,z)+w_\mu(z,y).\label{e:fla}
  \end{eqnarray}

A nonclassical example of ``generalized velocities'' satisfying \eq{e:axi}--\eq{e:axi3}
is given by: $w_\lambda(x,y)=\infty$ if $\lambda\le d(x,y)$, and
$w_\lambda(x,y)=0$ if $\lambda>d(x,y)$.

A ({\em metric\/}) {\em modular\/} on a set $X$ is any one-parameter family
$w=\{w_\lambda\}_{\lambda>0}$ of functions $w_\lambda:X\times X\to[0,\infty]$
satisfying \eq{e:axi}--\eq{e:axi3}. In particular, the family given by \eq{e:dxy} 
is the canonical (=\,natural) modular on a metric space $(X,d)$, which can be
interpreted as a field of average velocities.
For a different interpretation of modulars related to the joint generalized variation
of univariate maps and their relationships with classical modulars on linear spaces
we refer to~\cite{NA-I} (cf.\ also Section~\ref{s:exas}).

The difference between a metric (=\,distance function) and a modular on a set
is now clearly seen: a modular depends on a positive parameter and may assume
infinite values; the latter property means that it is impossible (or prohibited)
to move from $x$ to $y$ in time $\lambda$, unless one moves with infinite velocity
$w_\lambda(x,y)=\infty$. In addition (cf.\ \eq{e:dxy}), the ``velocity''
$w_\lambda(x,y)$ is {\em nonincreasing\/} as a function of ``time'' $\lambda>0$.  
The knowledge of ``average velocities'' $w_\lambda(x,y)$ for all $\lambda>0$ and
$x,y\in X$ provides more information than simply the knowledge of distances
$d(x,y)$ between $x$ and $y$: the distance $d(x,y)$ can be recovered as
a ``limit case'' via the formula (again cf.\ \eq{e:dxy}):
  $$d(x,y)=\inf\{\lambda>0:w_\lambda(x,y)\le1\}.$$

Now we describe briefly the main result of this paper. Given a modular $w$ on a set
$X$, we introduce the {\em modular space\/} $X_w^*=X_w^*(x_0)$ around a
point $x_0\in X$ as the set of those $x\in X$, for which $w_\lambda(x,x_0)$ is
finite for some $\lambda=\lambda(x)>0$. A map $T:X_w^*\to X_w^*$ is said
to be {\em modular contractive\/} if there exists a constant $0<k<1$ such that
for all small enough $\lambda>0$ and all $x,y\in X_w^*$ we have
$w_{k\lambda}(Tx,Ty)\le w_\lambda(x,y)$. Our main result (Theorem~\ref{t:main})
asserts that if $w$ is {\em convex\/} and {\em strict\/}, $X_w^*$ is
{\em modular complete\/} (the emphasized notions will be introduced in the main
text below) and $T:X_w^*\to X_w^*$ is modular contractive, then
$T$ admits a (unique) fixed point: $Tx_*=x_*$ for some $x_*\in X_w^*$.
The successive approximations of $x_*$ constructed in the proof of this result
converge to~$x_*$ in the modular sense, which is weaker than the metric convergence.
In particular, Banach's Contraction Principle follows if we take
into account~\eq{e:dxy}.

\smallbreak
This paper is organized as follows. In Section~\ref{s:mm} we study modulars
and convex modulars and introduce two modular spaces.
In Section~\ref{s:seq} we introduce the notions of modular convergence,
modular limit and modular completeness and show that they are ``weaker'' than
the corresponding metric notions. These notions are illustrated in Section~\ref{s:exas}
by examples. Section~\ref{s:fp} is devoted to a fixed point theorem for
modular contractions in modular complete modular metric spaces. This theorem
is then applied in Section~\ref{s:appl} to the existence of solutions of a
Carath{\'e}odory-type ordinary differential equation with the right-hand side from
the Orlicz space~$\mbox{L}^\vfi$. Finally, in Section~\ref{s:cr} some concluding
remarks are presented.

\section{Modulars and modular spaces} \label{s:mm}

In what follows $X$ is a nonempty set, $\la>0$ is understood in the
sense that $\la\in(0,\infty)$ and, in view of the disparity of the arguments,
functions $w:(0,\infty)\times X\times X\to[0,\infty]$ will be also written as
$w_\la(x,y)=w(\la,x,y)$ for all $\la>0$ and $x,\,y\in X$, so that
$w=\{w_\lambda\}_{\lambda>0}$ with $w_\lambda:X\times X\to[0,\infty]$.

\begin{defin}[\cite{DAN06,NA-I}] \label{d:md}
A function $w:(0,\infty)\times X\times X\to[0,\infty]$ is said to be a (metric)
{\em modular on $X$} if it satisfies the following three conditions:
\par(i) given $x,\,y\in X$, $x=y$ iff $w_\la(x,y)=0$ for all $\la>0$;
\par(ii) $w_\la(x,y)=w_\la(y,x)$ for all $\la>0$ and $x,\,y\in X$;
\par(iii) $w_{\la+\mu}(x,y)\le w_\la(x,z)+w_\mu(y,z)$ for all
$\la,\,\mu>0$ and $x,\,y,\,z\in X$.

\noindent If, instead of (i), the function $w$ satisfies only
\par(i${}'$) $w_\la(x,x)=0$ for all $\la>0$ and $x\in X$,\newline
then $w$ is said to be a {\em pseudomodular\/} on~$X$,
and if $w$ satisfies (i${}'$) and
\par(\is) given $x,y\in X$, if there exists a number $\lambda>0$,
possibly depending  \newline\hphantom{aaaaaa} on $x$ and $y$,
such that $w_\lambda(x,y)=0$, then $x=y$,\newline
the function $w$ is called a {\em strict modular\/} on~$X$.

A modular (pseudomodular, strict modular) $w$ on $X$ is said to be {\em convex\/} if,
instead if (iii), for all $\lambda,\mu>0$ and $x,y,z\in X$ it satisfies the inequality
\par\bigbreak
(iv) $\D w_{\la+\mu}(x,y)\le\frac\la{\la\!+\!\mu}\,w_\la(x,z)%
  +\frac\mu{\la\!+\!\mu}\,w_\mu(y,z)$.
\end{defin}

\medbreak
A motivation of the notion of {\em convexity} for modulars, which may look
unexpected at first glance, was given in \cite[Theorem~3.11]{NA-I}, cf.\ also
inequality \eq{e:fla}; a further generalization of this notion was presented
in~\cite[Section~5]{folia08}.

\smallbreak
Given a metric space $(X,d)$ with metric $d$, two {\em canonical\/} strict modulars are
associated with it: $w_\lambda(x,y)=d(x,y)$ (denoted simply by $d$),
which is independent of the first argument $\lambda$ and is a (nonconvex)
 modular on $X$ in the sense of (i)--(iii), and the {\em convex\/} modular \eq{e:dxy},
which satisfies (i), (ii) and (iv).
Both modulars $d$ and \eq{e:dxy} assume only finite values on~$X$.

\smallbreak
Clearly, if $w$ is a strict modular, then $w$ is a modular, which in turn implies
$w$ is a pseudomodular on~$X$, and similar implications hold for convex~$w$.

\smallbreak
The essential property of a pseudomodular $w$ on $X$ (cf.\ \cite[Section~2.3]{NA-I})
is that, for any given $x,y\in X$, the function
$0<\lambda\mapsto w_\lambda(x,y)\in[0,\infty]$ is {\em nonincreasing\/}
on $(0,\infty)$, and so, the limit from the right $w_{\lambda+0}(x,y)$ and
the limit from the left $w_{\lambda-0}(x,y)$ exist in $[0,\infty]$ and satisfy
the inequalities:
  \begin{equation} \label{e:lapm}
w_{\la+0}(x,y)\le w_\la(x,y)\le w_{\la-0}(x,y).
  \end{equation}

A {\em convex\/} pseudomodular $w$ on $X$ has the following additional property:
given $x,y\in X$, we have (cf.\ \cite[Section~3.5]{NA-I}):
  \begin{equation} \label{e:mpc}
\mbox{if \,$0<\mu\le\lambda$, \,then \,$w_\lambda(x,y)\le\D\frac\mu\lambda\,
w_\mu(x,y)\le w_\mu(x,y)$,}
  \end{equation}
i.e., functions $\lambda\mapsto w_\lambda(x,y)$ and
$\lambda\mapsto\lambda w_\lambda(x,y)$ are {\em nonincresing\/} on $(0,\infty)$.
%Moreover, by virtue of \eq{e:mpc},
%  $$%\begin{equation} \label{e:spine}
%\mbox{$0<\mu\le\lambda$ \,and \,$w_\mu(x,y)\le\lambda$ \,imply
% \,$w_\lambda(x,y)\le\mu$.}
%  $$%\end{equation}

\smallbreak
Throughout the paper we fix an element $x_0\in X$ arbitrarily.

\begin{defin}[\cite{DAN06,NA-I}] \label{d:modsp}
Given a pseudomodular $w$ on $X$, the two sets
  $$X_w\equiv X_w(x_0)=\bigl\{x\in X:
     \mbox{$w_\lambda(x,x_0)\to0$ \,as \,$\lambda\to\infty$}\bigr\}$$
\par\vspace{-10pt}\noindent 
and
\par\vspace{-20pt}
  $$X_w^*\equiv X_w^*(x_0)=\bigl\{x\in X:\mbox{$\exists\,\lambda=\lambda(x)>0$
     such that $w_\lambda(x,x_0)<\infty$}\bigr\}$$
are said to be {\em modular spaces\/} (around $x_0$).
\end{defin}

\smallbreak
It is clear that $X_w\subset X_w^*$, and it is known (cf.\  
\cite[Sections~3.1,\,3.2]{NA-I}) that this inclusion is proper in general.
It follows from \cite[Theorem~2.6]{NA-I} that if $w$ is a {\em modular\/} on $X$,
then the modular space $X_w$ can be equipped with a (nontrivial) metric
$d_w$, generated by $w$ and given by
  \begin{equation} \label{e:dw}
d_w(x,y)=\inf\{\lambda>0:w_\lambda(x,y)\le\lambda\},\quad\,\,x,\,y\in X_w.
  \end{equation}
It will be shown later that $d_w$ is a well defined metric on a larger set~$X_w^*$.

If $w$ is a {\em convex\/} modular on $X$, then according to
\cite[Section~3.5 and Theorem~3.6]{NA-I} the two modular spaces coincide,
$X_w=X_w^*$, and this common set can be endowed with a metric $d_w^*$ given by
  \begin{equation} \label{e:dw*}
d_w^*(x,y)=\inf\{\lambda>0:w_\lambda(x,y)\le1\},\quad\,\,x,\,y\in X_w^*;
  \end{equation}
moreover, $d_w^*$ is {\em specifically\/} equivalent to $d_w$
(see \cite[Theorem~3.9]{NA-I}).
By the convexity of $w$, the function
$\widehat w_\lambda(x,y)=\lambda w_\lambda(x,y)$ is a modular  on $X$
in the sense of (i)--(iii) and (cf.\ \cite[formula (3.3)]{NA-I})
  \begin{equation} \label{e:inclu}
X_{\widehat w}^*=X_w^*=X_w\supset X_{\widehat w},
  \end{equation}
where the last inclusion may be proper; moreover, $d_{\widehat w}=d_w^*$ on
$X_{\widehat w}$.

\smallbreak
Even if $w$ is a nonconvex modular on $X$, the quantity \eq{e:dw*} is also defined
for all $x,y\in X_w^*$, but it has only few properties
(cf.\ \cite[Theorem~3.6]{NA-I}): $d_w^*(x,x)=0$ and $d_w^*(x,y)=d_w^*(y,x)$.
In this case we have (cf.\ \cite[Theorem~3.9 and Example~3.10]{NA-I}):
if $d_w(x,y)<1$, then $d_w^*(x,y)\le d_w(x,y)$,
and if $d_w^*(x,y)\ge1$, then $d_w(x,y)\le d_w^*(x,y)$.

\smallbreak
Let us illustrate the above in the case of a metric space $(X,d)$ with the two
canonical modulars $d$ and $w$ from \eq{e:dxy} on it. We have:
$X_d=\{x_0\}\subset X_d^*=X_w=X_w^*=X$, and given $x,y\in X$,
$d_d(x,y)=d(x,y)$, $d_d^*(x,y)=0$, $d_w(x,y)=\sqrt{d(x,y)}$,
$d_w^*(x,y)=d(x,y)$ and $\widehat d(x,y)=\lambda w_\lambda(x,y)=d(x,y)$.
Thus, the convex modular $w$ from \eq{e:dxy} plays a more adequate role
in restoring the metric space $(X,d)$ from $w$ (cf.\ $d_w^*=d$ on
$X_w=X_w^*=X$, whereas $X_d\subset X_d^*=X$, $d_d=d$ and $d_d^*=0$),
and so, in what follows any metric space $(X,d)$ will be considered equipped
only with the modular~\eq{e:dxy}. This convention is also justified as follows.

\smallbreak
Now we exhibit the relationship between convex and nonconvex modulars and
show that $d_w$ is a well defined metric on $X_w^*$ (and not only on~$X_w$).
If $w$ is a (not necessarily convex) modular on $X$, then the function
(cf.\ \eq{e:dxy} where $d(x,y)$ plays the role of a modular)
  $$v_\lambda(x,y)=\frac{w_\lambda(x,y)}{\lambda},\qquad\lambda>0,
     \quad x,y\in X,$$
is always a {\em convex\/} modular on $X$. In fact, conditions (i) and (ii) are
clear for $v$ and, as for (iv), we have, by virtue of (iii) for $w$:
  \begin{eqnarray}
v_{\lambda+\mu}(x,y)\!\!\!&=\!\!\!&\D\frac{w_{\lambda+\mu}(x,y)}{\lambda+\mu}\le
  \frac{w_\lambda(x,z)+w_\mu(y,z)}{\lambda+\mu}=\nonumber\\
\!\!\!&=\!\!\!&\frac{\lambda}{\lambda\!+\!\mu}\!\cdot\!\frac{w_\lambda(x,z)}{\lambda}\!+\!
  \frac{\mu}{\lambda\!+\!\mu}\!\cdot\!\frac{w_\mu(y,z)}{\mu}=
  \frac{\lambda}{\lambda\!+\!\mu}\,v_\lambda(x,z)\!+\!
  \frac{\mu}{\lambda\!+\!\mu}\,v_\mu(y,z).\nonumber
  \end{eqnarray}
Moreover, because $w=\widehat v$, we find from \eq{e:inclu} that
$X_w\subset X_w^*=X_v=X_v^*$. Since
$d_v^*(x,y)=\inf\{\lambda>0:w_\lambda(x,y)/\lambda\le1\}=d_w(x,y)$
for all $x,y\in X_w^*$, i.e., $d_v^*=d_w$ on $X_w^*$, and $d_v^*$ is a
metric on $X_v^*=X_w^*$, then we conclude that $d_w$ is a
{\em well defined metric on $X_w^*$} (the same conclusion follows immediately
from \cite[Theorem~1]{folia08}) with $X'=X_w^*$). This property distinguishes
our theory of modulars from the classical theory: if $\rho$ is a classical modular
on a linear space $X$ in the sense of Musielak and Orlicz (\cite{Musielak}) and
$w_\lambda(x,y)\!=\!\rho((x-y)/\lambda)$, \mbox{$\lambda\!>\!0$}, $x,y\in X$,
 then the expression
$v_\lambda(x,y)=(1/\lambda)w_\lambda(x,y)=(1/\lambda)\rho((x-y)/\lambda)$
is {\em not allowed\/} as a classical modular on~$X$.
Since $v$ is convex and $d_v^*=d_w$ on $X_w^*$, given $x,y\in X_w^*$,
by virtue of \cite[Theorem~3.9]{NA-I}, we have:
  $$d_w(x,y)<1\,\,\,\mbox{iff}\,\,\,d_v(x,y)<1,\,\mbox{and}\,\,\,
     d_w(x,y)\le d_v(x,y)\le\sqrt{d_w(x,y)};$$
  $$d_w(x,y)\ge1\,\,\,\mbox{iff}\,\,\,d_v(x,y)\ge1,\,\mbox{and}\,\,\,
     \sqrt{d_w(x,y)}\le d_v(x,y)\le d_w(x,y).$$

More metrics can be defined on $X_w^*$ for a given modular $w$ on $X$
in the following general way (cf.\ \cite[Theorem~1]{folia08}): if $\RB^+=[0,\infty)$
and $\kappa:\RB^+\to\RB^+$ is superadditive (i.e., 
$\kappa(\lambda)+\kappa(\mu)\le\kappa(\lambda+\mu)$ for all
$\lambda,\mu\ge0$) and such that $\kappa(u)>0$ for $u>0$ and
$\kappa(+0)=\lim_{u\to+0}\kappa(u)=0$, then the function
$d_{\kappa,w}(x,y)=\inf\{\lambda>0:w_\lambda(x,y)\le\kappa(\lambda)\}$
is a well defined metric on~$X_w^*$.
 
\smallbreak
Given a pseudomodular (modular, strict modular, convex or not) $w$ on $X$,
$\lambda>0$ and $x,y\in X$, we define the {\em left\/} and {\em right\/}
{\em regularizations} of $w$ by
  $$w_\lambda^-(x,y)=w_{\lambda-0}(x,y)\quad\,\,\mbox{and}
     \quad\,\,w_\lambda^+(x,y)=w_{\lambda+0}(x,y).$$
Since, by \eq{e:lapm}, $w_\lambda^+(x,y)\le w_\lambda(x,y)\le w_\lambda^-(x,y)$, and
  \begin{equation} \label{e:la2la1}
w_{\lambda_2}^-(x,y)\le w_\lambda(x,y)\le w_{\lambda_1}^+(x,y)
     \quad\mbox{for \,all}\quad 0<\lambda_1<\lambda<\lambda_2,
  \end{equation}
it is a routine matter to verify that $w^-$ and $w^+$ are pseudomodulars
(modulars, strict modulars, convex or not, respectively) on $X$,
$X_{w^-}\!=X_w=X_{w^+}$, $X_{w^-}^*\!=X_w^*=X_{w^+}^*$,
$d_{w^-}\!=d_w=d_{w^+}$ on $X_w$ and
$d_{w^-}^*\!=d_w^*=d_{w^+}^*$ on~$X_w^*$.
For instance, let us check the last two equalities for metrics.
Given $x,y\in X_w^*$, by virtue of \eq{e:lapm}, we find
\mbox{$d_{w^-}^*(x,y)\ge d_w^*(x,y)\ge d_{w^+}^*(x,y)$}.
In order to see that $d_{w^-}^*(x,y)\le d_w^*(x,y)$, we let $\lambda>d_w^*(x,y)$
be arbitrary, choose $\mu$ such that $d_w^*(x,y)<\mu<\lambda$, which, by
\eq{e:la2la1}, gives $w_\lambda^-(x,y)\le w_\mu(x,y)\le1$, and so,
$d_{w^-}^*(x,y)\le\lambda$, and then let $\lambda\to d_w^*(x,y)$.
In order to prove that $d_w^*(x,y)\le d_{w^+}^*(x,y)$, we let
$\lambda>d_{w^+}^*(x,y)$ be arbitrary, choose $\mu$ such that
$d_{w^+}^*(x,y)<\mu<\lambda$, which, by \eq{e:la2la1}, implies
$w_\lambda(x,y)\le w_\mu^+(x,y)\le1$, and so,
$d_w^*(x,y)\le\lambda$, and then let $\lambda\to d_{w^+}^*(x,y)$.

\smallbreak
In this way we have seen that the regularizations provide no new modular spaces
as compared to $X_w$ and $X_w^*$ and no new metrics as compared to $d_w$
and $d_w^*$. The right regularization will be needed in Section~\ref{s:fp}
for the characterization of metric Lipschitz maps in terms of underlying modulars.

\section{Sequences in modular spaces and modular convergence} \label{s:seq}

The notions of modular convergence, modular limit, modular completeness, etc.,
which we study in this section, are known in the classical theory of modulars on
linear spaces (e.g., \cite{Maligranda,Musielak,Orlicz,RaoRen}).  Since the theory
of (metric) modulars from \cite{DAN06}--\cite{NA-II} is significantly more general
than the classical theory, the notions mentioned above do not carry over to metric
modulars in a straightforward way and ought to be reintroduced and justified.

\begin{defin} \label{d:wconv}
Given a pseudomodular $w$ on $X$, a sequence
of elements $\{x_n\}\equiv\{x_n\}_{n=1}^\infty$ from $X_w$ or $X_w^*$
is said to be {\em modular convergent\/} (more precisely, {\em $w$-convergent\/}) to an
element $x\in X$ if there exists a number $\lambda>0$, possibly depending on
$\{x_n\}$ and $x$, such that $\lim_{n\to\infty}w_\lambda(x_n,x)=0$.
This will be written briefly as $x_n\wto x$ (as $n\to\infty$), and any such
element $x$ will be called a {\em modular limit\/} of the sequence~$\{x_n\}$.
\end{defin}

Note that if $\lim_{n\to\infty}w_\lambda(x_n,x)=0$, then by virtue
of the monotonicity of the function $\lambda'\mapsto w_{\lambda'}(x_n,x)$,
we have: $\lim_{n\to\infty}w_\mu(x_n,x)=0$ for all $\mu\ge\lambda$.

\smallbreak
It is clear for a metric space $(X,d)$ and the modular \eq{e:dxy} on it that
the metric convergence and the modular convergence in $X$ coincide.

\smallbreak
We are going to show that the modular convergence is much weaker than the
metric convergence (in the sense to be made more precise below).
First, we study to what extent the above definition is correct, and what is the
relationship between the modular and metric convergences in $X_w$ and $X_w^*$.

\begin{thm} \label{t:corr}
Let $w$ be a pseudomodular on $X$. We have\/{\rm:}

{\rm(a)} the modular spaces $X_w$ and $X_w^*$ are closed with respect
to the modular convergence, i.e., if\/ $\{x_n\}\subset X_w$ {\rm(}or $X_w^*${\rm)},
$x\in X$ and $x_n\wto x$, then $x\in X_w$ {\rm(}or $x\in X_w^*$, respectively\/{\rm);}

{\rm(b)} if\/ $w$ is a strict modular on $X$, then the modular limit
is determined uniquely {\rm(}if it exists\/{\rm)}.
\end{thm}
\begin{pf}
(a) Since $x_n\wto x$, there exists a $\lambda_0=\lambda_0(\{x_n\},x)>0$
such that $w_{\lambda_0}(x_n,x)\to0$ as $n\to\infty$.

1. First we treat the case when $\{x_n\}\subset X_w$. Let $\vep>0$ be arbitrarily
fixed. Then there is an $n_0=n_0(\vep)\in\NB$ such that
$w_{\lambda_0}(x_{n_0},x)\le\vep/2$. Since $x_{n_0}\in X_w=X_w(x_0)$,
we have $w_\lambda(x_{n_0},x_0)\to0$ as $\lambda\to\infty$, and so,
there exists a $\lambda_1=\lambda_1(\vep)>0$ such that
$w_{\lambda_1}(x_{n_0},x_0)\le\vep/2$. Then conditions (iii) and (ii) from
Definition~\ref{d:md} imply
  $$w_{\lambda_0+\lambda_1}(x,x_0)\le w_{\lambda_0}(x,x_{n_0})+
     w_{\lambda_1}(x_0,x_{n_0})\le\vep.$$
The function $\lambda\mapsto w_\lambda(x,x_0)$ is nonincreasing on $(0,\infty)$,
and so,
  $$w_\lambda(x,x_0)\le w_{\lambda_0+\lambda_1}(x,x_0)\le\vep\quad
     \mbox{for \,all}\quad\lambda\ge\lambda_0+\lambda_1,$$
implying $w_\lambda(x,x_0)\to0$ as $\lambda\to\infty$, i.e., $x\in X_w$.

2. Now suppose that $\{x_n\}\subset X_w^*$. Then there exists an $n_0\in\NB$
such that $w_{\lambda_0}(x_{n_0},x)\le1$. Since $x_{n_0}\in X_w^*=X_w^*(x_0)$,
there is a $\lambda_1>0$ such that $w_{\lambda_1}(x_{n_0},x_0)<\infty$.
Now it follows from conditions (iii) and (ii) that
  $$w_{\lambda_0+\lambda_1}(x,x_0)\le w_{\lambda_0}(x,x_{n_0})+
     w_{\lambda_1}(x_0,x_{n_0})<\infty,$$
and so, $x\in X_w^*$.

\smallbreak
(b) Let $\{x_n\}\subset X_w$ or $X_w^*$ and $x,y\in X$ be such that
$x_n\wto x$ and $x_n\wto y$. By the definition of the modular convergence,
there exist $\lambda=\lambda(\{x_n\},x)>0$ and $\mu=\mu(\{x_n\},y)>0$
such that $w_\lambda(x_n,x)\to0$ and $w_\mu(x_n,y)\to0$ as $n\to\infty$.
By conditions (iii) and (ii),
  $$w_{\lambda+\mu}(x,y)\le w_\lambda(x,x_n)+w_\mu(y,x_n)\to0\quad
    \mbox{as}\quad n\to\infty.$$
It follows that $w_{\lambda+\mu}(x,y)=0$, and so, by condition (\is) from
Definition~\ref{d:md}, we get $x=y$.
\qed\end{pf}

\medbreak

It was shown in \cite[Theorem~2.13]{NA-I} that if $w$ is a modular on $X$,
then for $\{x_n\}\subset X_w$ and $x\in X_w$ we have:
  \begin{equation} \label{e:dwwla}
\lim_{n\to\infty}d_w(x_n,x)=0\quad\,\mbox{iff}\quad\,
\lim_{n\to\infty}w_\lambda(x_n,x)=0\,\,\mbox{for \,all}\,\,\lambda>0.
  \end{equation}
and so, the metric convergence (with respect to the metric $d_w$) implies the modular
convergence (cf.\ Definition~\ref{d:wconv}), but not vice versa in general.
As the proof of \cite[Theorem~2.13]{NA-I} suggests, \eq{e:dwwla} is also true for
$\{x_n\}\subset X_w^*$ and $x\in X_w^*$. An assertion similar to \eq{e:dwwla} holds for Cauchy sequences from the modular spaces $X_w$ and $X_w^*$.

\smallbreak
Now we establish a result similar to \eq{e:dwwla} for {\em convex\/} modulars.

\begin{thm} \label{t:dil=w}
Let $w$ be a convex modular on $X$. Given a sequence $\{x_n\}$ from
$X_w^*\,(\,=X_w)$ and an element $x\in X_w^*$, we have\/{\rm:}
  $$\lim_{n\to\infty}d_w^*(x_n,x)=0\quad\,\mbox{iff}\quad\,
     \lim_{n\to\infty}w_\lambda(x_n,x)=0\,\,\,\mbox{for all}\,\,\,\lambda>0.$$
\par A similar assertion holds for Cauchy sequences with respect to~$d_w^*$.
\end{thm}
\begin{pf}
Step~1. {\em Sufficiency.} Given $\vep>0$, by the assumption, there exists
a number $n_0(\vep)\in\NB$ such that $w_\vep(x_n,x)\le1$ for all $n\ge n_0(\vep)$,
and so, the definition \eq{e:dw*} of $d_w^*$ implies $d_w^*(x_n,x)\le\vep$
for all $n\ge n_0(\vep)$.

\smallbreak
{\em Necessity.} First, suppose that $0<\lambda\le1$. Given $\vep>0$, we have:
either (a) $\vep<\lambda$, or (b) $\vep\ge\lambda$. In case (a), by the
assumption, there is an $n_0(\vep)\in\NB$ such that $d_w^*(x_n,x)<\vep^2$
for all $n\ge n_0(\vep)$, and so, by the definition of $d_w^*$,
$w_{\vep^2}(x_n,x)\le1$ for all $n\ge n_0(\vep)$. Since $\vep^2<\lambda^2\le\lambda$
and $\vep<\lambda$, inequality \eq{e:mpc} yields:
  $$w_\lambda(x_n,x)\le\frac{\vep^2}{\lambda}\,w_{\vep^2}(x_n,x)\le
     \frac{\vep}{\lambda}\,\vep<\vep\quad\mbox{for all}\quad n\ge n_0(\vep).$$
In case (b) we set $n_1(\vep)=n_0(\lambda/2)$, where $n_0(\cdot)$ is as above.
Then, as we have just established,
$w_\lambda(x_n,x)<\lambda/2\le\vep/2<\vep$ for all $n\ge n_1(\vep)$.

\smallbreak
Now, assume that $\lambda>1$. Again, given $\vep>0$, we have:
either (a) $\vep<\lambda$, or (b) $\vep\ge\lambda$. In case (a) there is an
$N_0(\vep)\in\NB$ such that $d_w^*(x_n,x)<\vep$ for all $n\ge N_0(\vep)$,
and so, $w_\vep(x_n,x)\le1$ for all $n\ge N_0(\vep)$. Since $\vep<\lambda$
and $\lambda>1$, by virtue of \eq{e:mpc}, we find
  $$w_\lambda(x_n,x)\le\frac{\vep}{\lambda}\,w_\vep(x_n,x)\le
     \frac\vep\lambda<\vep\quad\mbox{for all}\quad n\ge N_0(\vep).$$
In case (b) we put $N_1(\vep)=N_0(\lambda/2)$, where $N_0(\cdot)$ is as above.
Then it follows that $w_\lambda(x_n,x)<\lambda/2\le\vep/2<\vep$ for all
$n\ge N_1(\vep)$.

Thus, we have shown that $w_\lambda(x_n,x)\to0$ as $n\to\infty$ for all $\lambda>0$.

\smallbreak
Step~2. The assertion for Cauchy sequences is of the form:
  $$\lim_{n,m\to\infty}d_w^*(x_n,x_m)=0\quad\,\mbox{iff}\quad\,
     \lim_{n,m\to\infty}w_\lambda(x_n,x_m)=0\,\,\,\mbox{for all}\,\,\,\lambda>0;$$
its proof is similar to the one given in Step~1 with suitable modifications.
\qed\end{pf}

Theorem~\ref{t:dil=w} shows, in particular, that in a metric space $(X,d)$ with
modular \eq{e:dxy} on it the metric and modular convergences are equivalent.

\begin{defin} \label{d:Delta2}
A pseudomodular $w$ on $X$ is said to {\em satisfy the\/} (sequential)
{\em $\Delta_2$-condition\/} (on $X_w^*$) if the following condition holds:
given a sequence $\{x_n\}\subset X_w^*$ and $x\in X_w^*$, if there exists
a number $\lambda>0$, possibly depending on $\{x_n\}$ and $x$, such that
$\lim_{n\to\infty}w_\lambda(x_n,x)=0$, then
$\lim_{n\to\infty}w_{\lambda/2}(x_n,x)=0$.
\end{defin}

A similar definition applies with $X_w^*$ replaced by $X_w$.

\smallbreak
In the case of a metric space $(X,d)$ the modular \eq{e:dxy} clearly satisfies
the $\Delta_2$-condition on~$X$.

\smallbreak
The following important observation, which generalizes the corresponding result
from the theory of classical modulars on linear spaces (cf.\ \cite[I,5.2.IV]{Musielak}),
provides a criterion for the metric and modular convergences to coincide.

\begin{thm} \label{t:dew}
Given a modular $w$ on $X$, we have\/{\rm:} the metric convergence on $X_w^*$
{\rm(}with respect to $d_w$ if $w$ is arbitrary, and with respect to $d_w^*$
if $w$ is convex\/{\rm)} coincides with the modular convergence
iff $w$ satisfies the \mbox{$\Delta_2$-condition}~on~$X_w^*$.
\end{thm}
\begin{pf}
Let $\{x_n\}\subset X_w^*$ and $x\in X_w^*$ be given. We know from \eq{e:dwwla}
and Theorem\/~\ref{t:dil=w} that the metric convergence (with respect to $d_w$
if $w$ is a modular or with respect to $d_w^*$ if $w$ is a convex modular) of
$x_n$ to $x$ is equivalent to
  \begin{equation} \label{e:woc}
\lim_{n\to\infty}w_\lambda(x_n,x)=0\quad\,\mbox{for \,all}\quad\,\lambda>0.
  \end{equation}

($\Rightarrow\!$) Suppose that the metric convergence coincides with the modular
convergence on $X_w^*$. If there exists a $\lambda_0>0$ such that
$w_{\lambda_0}(x_n,x)\to0$ as $n\to\infty$, then $x_n$ is modular convergent
to~$x$, and so, $x_n$ converges to $x$ in metric ($d_w$ or $d_w^*$).
It follows that \eq{e:woc} holds implying, in particular,
 $w_{\lambda_0/2}(x_n,x)\to0$ as
$n\to\infty$, and so, $w$ satisfies the $\Delta_2$-condition.

\smallbreak
($\Leftarrow$) By virtue of \eq{e:woc}, the metric convergence on $X_w^*$
always implies the modular convergence, and so, it suffices to verify the converse
assertion, namely: if $x_n\wto x$, then \eq{e:woc} holds. In fact, if
$x_n\wto x$, then $w_{\lambda_0}(x_n,x)\to0$ as $n\to\infty$ for some
constant $\lambda_0=\lambda_0(\{x_n\},x)>0$. The $\Delta_2$-condition implies
$w_{\lambda_0/2}(x_n,x)\to0$ as $n\to\infty$, and so, the induction yields
$w_{\lambda_0/2^j}(x_n,x)\to0$ as $n\to\infty$ for all $j\in\NB$. Now,
given $\lambda>0$, there exists a $j=j(\lambda)\in\NB$ such that
$\lambda>\lambda_0/2^j$. By the monotonicity of $\lambda\mapsto w_\lambda(x_n,x)$,
we have:
  $$w_\lambda(x_n,x)\le w_{\lambda_0/2^j}(x_n,x)\to0\quad
     \mbox{as}\quad n\to\infty.$$
By the arbitrariness of $\lambda>0$, condition \eq{e:woc} follows.
\qed\end{pf}

\begin{defin} \label{d:wCauchy}
Given a modular $w$ on $X$, a sequence $\{x_n\}\subset X_w^*$ is said to be
{\em modular Cauchy\/} (or {\em $w$-Cauchy\/}) if there exists a number
$\lambda=\lambda(\{x_n\})>0$ such that $w_\lambda(x_n,x_m)\to0$
as $n,m\to\infty$, i.e.,
  $$\mbox{$\forall\,\vep>0$ $\exists\,n_0(\vep)\in\NB$ such that 
     $\forall\,n\ge n_0(\vep)$, $m\ge n_0(\vep)$: $w_\lambda(x_n,x_m)\le\vep$.}$$
\end{defin}

It follows from Theorem~\ref{t:dil=w} (Step~2 in its proof) and
Definition~\ref{d:wCauchy} that a sequence from $X_w^*$, which is
Cauchy in metric $d_w$ or $d_w^*$, is modular Cauchy.

\smallbreak
Note that a modular convergent sequence is modular Cauchy. In fact, if
$x_n\wto x$, then $w_\lambda(x_n,x)\to0$ as $n\to\infty$ for some $\lambda>0$,
and so, for each $\vep>0$ there exists an $n_0(\vep)\in\NB$ such that
$w_\lambda(x_n,x)\le\vep/2$ for all $n\ge n_0(\vep)$. It follows from (iii) that
if $n,m\ge n_0(\vep)$, then
$w_{2\lambda}(x_n,x_m)\le w_\lambda(x_n,x)+w_\lambda(x_m,x)\le\vep$,
which implies that $\{x_n\}$ is modular Cauchy.

\smallbreak
The following definition will play an important role below.
\begin{defin} \label{d:wcomp}
Given a modular $w$ on $X$, the modular space $X_w^*$ is said to be
{\em modular complete\/} (or {\em $w$-complete\/}) if each modular Cauchy
sequence from $X_w^*$ is modular convergent in the following (more precise) sense:
if $\{x_n\}\subset X_w^*$ and there exists a $\lambda=\lambda(\{x_n\})>0$
such that $\lim_{n,m\to\infty}w_\lambda(x_n,x_m)=0$, then there exists
an $x\in X_w^*$ such that $\lim_{n\to\infty}w_\lambda(x_n,x)=0$.
\end{defin}

The notions of modular convergence, modular limit and modular completeness,
introduced above, are illustrated by examples in the next section. It is clear
from \eq{e:dxy} that for a metric space $(X,d)$ these notions coincide with
respective notions in the metric space setting.

\section{Examples of metric and modular convergences} \label{s:exas}

We begin with recalling certain properties of $\vfi$-functions and convex functions
on the set of all nonnegative reals $\RB^+=[0,\infty)$.

\smallbreak
A function $\vfi:\RB^+\to\RB^+$ is said to be a {\em $\vfi$-function\/} if
it is continuous, nondecreasing, unbounded (and so,
$\vfi(\infty)\equiv\lim_{u\to\infty}\vfi(u)=\infty$) and assumes
the value zero only at zero: $\vfi(u)=0$ iff $u=0$. 

\smallbreak
If $\vfi:\RB^+\to\RB^+$ is a convex function such that $\vfi(u)=0$ iff $u=0$, then
it is (automatically) continuous, strictly increasing and unbounded, and so, it is a
convex $\vfi$-function. Also, $\vfi$ is superadditive:
$\vfi(u_1)+\vfi(u_2)\le\vfi(u_1+u_2)$ for all $u_1,u_2\in\RB^+$
(cf.\ \cite[Section~I.1]{KR}). Moreover, $\vfi$ admits the inverse function
$\vfi^{-1}:\RB^+\to\RB^+$, which is continuous, strictly increasing,
$\vfi^{-1}(u)=0$ iff $u=0$, $\vfi^{-1}(\infty)=\infty$, and which is subadditive:
$\vfi^{-1}(u_1+u_2)\le\vfi^{-1}(u_1)+\vfi^{-1}(u_2)$
for all $u_1,u_2\in\RB^+$. The function $\vfi$ is said to {\em satisfy the
$\Delta_2$-condition at infinity\/} (cf.\ \cite[Section~I.4]{KR}) if there exist
constants $K>0$ and $u_0\ge0$ such that $\vfi(2u)\le K\vfi(u)$ for all $u\ge u_0$.

\smallbreak
{\bf\ref{s:exas}.1.} Let the triple $(M,d,+)$ be a {\em metric semigroup}, i.e.,
the pair $(M,d)$ is a metric space with metric $d$, the pair $(M,+)$ is an Abelian
semigroup with respect to the operation of addition $+$, and $d$ is translation
invariant in the sense that $d(p+r,q+r)=d(p,q)$ for all $p,q,r\in M$. Any normed
linear space $(M,|\cdot|)$ is a metric semigroup with the induced metric
$d(p,q)=|p-q|$, $p,q\in M$, and the addition operation $+$ from $M$.
If $K\subset M$ is a convex cone (i.e., $p+q,\lambda p\in K$ whenever $p,q\in K$
and $\lambda\ge0$), then the triple $(K,d,+)$ is also a metric semigroup.
A nontrivial example of a metric semigroup is as follows (cf.\ \cite{DeBlasi,Rad}).
Let $(Y,|\cdot|)$ be a real normed space and $M$ be the family of all nonempty
closed bounded convex subsets of $Y$ equipped with the Hausdorff metric $d$
given by $d(P,Q)\!=\!\max\{\mbox{e}(P,Q),\mbox{e}(Q,P)\}$, where
\mbox{$P,Q\!\in\! M$} and $\mbox{e}(P,Q)=\sup_{p\in P}\inf_{q\in Q}|p-q|$.
Given $P,Q\in M$, we define $P\oplus Q$ as the closure in $Y$ of the Minkowski sum
$P+Q=\{p+q:\mbox{$p\in P$,\,$q\in Q$}\}$. Then the triple $(M,d,\oplus)$ is a
metric semigroup (actually, $M$ is an abstract convex cone).
For more information on metric semigroups and their special cases, abstract
convex cones, including examples we refer to
\cite{Sovae,folia04,NA-I,NA-II} and references therein.

Given a closed interval $[a,b]\subset\RB$ with $a<b$, we denote by $\XB=M^{[a,b]}$
the set of all mappings $x:[a,b]\to M$. If $\vfi$ is a {\em convex\/} $\vfi$-function
on $\RB^+$, we define a function $w:(0,\infty)\times\XB\times\XB\to[0,\infty]$ for
all $\lambda>0$ and $x,y\in\XB$ by (note that $w$ depends on~$\vfi$)
  \begin{equation} \label{e:wlaxy}
w_\lambda(x,y)=\sup_\pi\sum_{i=1}^m\vfi\Biggl(
\frac{d\bigl(x(t_i)+y(t_{i-1}),x(t_{i-1})+y(t_i)\bigr)}{\lambda\!\cdot\!(t_i-t_{i-1})}
\Biggr)\!\!\cdot\!(t_i-t_{i-1}),
  \end{equation}
where the supremum is taken over all partitions $\pi=\{t_i\}_{i=1}^m$ of the interval
$[a,b]$, i.e., $m\in\NB$ and $a=t_0<t_1<t_2<\dots<t_{m-1}<t_m=b$.
It was shown in \cite[Sections~3, 4]{Sovae} that $w$ is a {\em convex pseudomodular\/}
on $\XB$. Thus, given $x_0\in M$, the modular space $\XB_w^*=\XB_w^*(x_0)$
(here $x_0$ denotes also the constant mapping $x_0(t)=x_0$ for all $t\in[a,b]$),
which was denoted in \cite[(3.20) and Section~4.1]{Sovae} by $\GV\vfi([a,b];M)$
and called the {\em space of mappings of bounded generalized $\vfi$-variation},
is well defined and, by the translation invariance of $d$ on $M$, we have:
$x\in\XB_w^*=\GV\vfi([a,b];M)$ iff $x:[a,b]\to M$ and there exists a constant
$\lambda=\lambda(x)>0$ such that
  \begin{equation} \label{e:wlax0}
w_\lambda(x,x_0)=\sup_\pi\sum_{i=1}^m\vfi\Biggl(
\frac{d\bigl(x(t_i),x(t_{i-1})\bigr)}{\lambda\,(t_i-t_{i-1})}\Biggr)(t_i-t_{i-1})<\infty.
  \end{equation}
Note that $w_\lambda(x,x_0)$ from \eq{e:wlax0} is independent of $x_0\in M$; this
value is called the {\em generalized $\vfi_\lambda$-variation\/} of $x$, where
$\vfi_\lambda(u)=\vfi(u/\lambda)$, $u\in\RB^+$. Since $w$ satisfies on $\XB$
conditions (i${}'$), (ii) and (iv) (and not (i) in general) from definition~\ref{d:md},
the quantity $d_w^*$ from \eq{e:dw*} is only a pseudometric on $\XB_w^*$
and, in particular, only $d_w^*(x,x)=0$ holds for $x\in\XB_w^*$ (note that
$d_w^*(x,y)$ was denoted by $\Delta_\vfi(x,y)$ in \cite[equality (4.5)]{Sovae}).

\smallbreak
{\bf\ref{s:exas}.2.} In order to ``turn'' \eq{e:wlaxy} into a modular, we fix an
$x_0\in M$ and set $X=\{x:[a,b]\to M\mid x(a)=x_0\}\subset\XB$. We assert
that $w$ from \eq{e:wlaxy} is a {\em strict\/} convex modular on~$X$. In fact,
given $x,y\in X$ and $t,s\in[a,b]$ with $t\ne s$, it follows from \eq{e:wlaxy} that
  $$\vfi\Biggl(\frac{d\bigl(x(t)+y(s),x(s)+y(t)\bigr)}{\lambda\,|t-s|)}\Biggr)|t-s|
     \le w_\lambda(x,y),$$
and so, by the translation invariance of $d$ and the triangle inequality,
  \begin{eqnarray}
|d(x(t),y(t))-d(x(s),y(s))|\!\!&\le\!\!&d(x(t)\!+\!y(s),x(s)\!+\!y(t))\le\nonumber\\[2pt]
\!\!&\le\!\!&\lambda\,|t-s|\,\vfi^{-1}\Biggl(\frac{w_\lambda(x,y)}{|t-s|}\Biggr).
  \label{e:fwL}
  \end{eqnarray}
Now, if we suppose that $w_\lambda(x,y)=0$ for some $\lambda>0$, then
for all $t\in[a,b]$, $t\ne s=a$, we get (note that $x(a)=y(a)=x_0$)
  $$d(x(t),y(t))=|d(x(t),y(t))-d(x(a),y(a))|\le0.$$
Thus, $x(t)=y(t)$ for all $t\in[a,b]$, and so, $x=y$ as elements of~$X$.

\smallbreak
It is clear for the modular space $X_w^*=X_w^*(x_0)$ that
  \begin{equation} \label{e:GVX}
X_w^*=\XB_w^*\cap X=\GV\vfi([a,b];M)\cap X,
  \end{equation}
i.e., $x\in X_w^*$ iff $x:[a,b]\to M$, $x(a)=x_0$ and \eq{e:wlax0} holds
for some $\lambda>0$. Moreover, the function $d_w^*$ from \eq{e:dw*}
is a {\em metric\/} on~$X_w^*$.

\smallbreak
{\bf\ref{s:exas}.3.} In this subsection we show that if $(M,d,+)$ is a {\em complete\/}
metric semigroup (i.e., $(M,d)$ is complete as a metric space), then the modular space
$X_w^*$ from \eq{e:GVX} is {\em modular complete\/} in the sense
of definition~\ref{d:wcomp}.

\smallbreak
Let $\{x_n\}\subset X_w^*$ be a $w$-Cauchy sequence, so that
$w_\lambda(x_n,x_m)\to0$ as $n,m\to\infty$ for some constant
$\lambda=\lambda(\{x_n\})>0$. Given $n,m\in\NB$ and $t\in[a,b]$, $t\ne a$,
it follows from \eq{e:fwL} with $x=x_n$, $y=x_m$ and $s=a$ that
(again note that $x_n(a)=x_0$ for all $n\in\NB$)
  $$d(x_n(t),x_m(t))\le\lambda\,(t-a)\,\vfi^{-1}\Biggl(
     \frac{w_\lambda(x_n,x_m)}{t-a}\Biggr).$$
This estimate, the modular Cauchy property of~$\{x_n\}$, the continuity of~$\vfi^{-1}$
and the completeness of $(M,d,+)$ imply the existence of an $x:[a,b]\to M$,
$x(a)=x_0$ (and so, $x\in X$), such that the sequence $\{x_n\}$ converges pointwise
on $[a,b]$ to $x$, i.e., $\lim_{n\to\infty}d(x_n(t),x(t))=0$ for all $t\in[a,b]$.
We assert that $\lim_{n\to\infty}w_\lambda(x_n,x)=0$. By the (sequential) lower
semicontinuity of the functional $w_\lambda(\cdot,\cdot)$ from \eq{e:wlaxy}
(cf.\ \cite[assertion (4.8) on p.~27]{Sovae}), we get
  \begin{equation} \label{e:liminf}
w_\lambda(x_n,x)\le\liminf_{m\to\infty}w_\lambda(x_n,x_m)\quad
\mbox{for \,all}\quad n\in\NB.
  \end{equation}
Now, given $\vep>0$, by the modular Cauchy condition for $\{x_n\}$, there is
an $n_0(\vep)\in\NB$ such that $w_\lambda(x_n,x_m)\le\vep$ for all
$n\ge n_0(\vep)$ and $m\ge n_0(\vep)$, and so,
  $$\limsup_{m\to\infty}w_\lambda(x_n,x_m)\le
     \sup_{m\ge n_0(\vep)}w_\lambda(x_n,x_m)\le\vep\quad
     \mbox{for \,all}\quad n\ge n_0(\vep).$$
Since the limit inferior does not exceed the limit superior (for any real sequences),
it follows from the last displayed line and \eq{e:liminf} that $w_\lambda(x_n,x)\le\vep$
for all $n\ge n_0(\vep)$, i.e., $w_\lambda(x_n,x)\to0$ as $n\to\infty$. Finally,
since, by Theorem~\ref{t:corr}\,(a), $X_w^*$  is closed with respect to the modular
convergence, we infer that $x\in X_w^*$, which was to be proved.

\smallbreak
{\bf\ref{s:exas}.4.} In order to be able to calculate explicitly, for the sake of simplicity
we assume furthermore that $M=\RB$ with $d(p,q)=|p-q|$, $p,q\in\RB$, and the
function $\vfi$ satisfies the {\em Orlicz condition at infinity\/}: $\vfi(u)/u\to\infty$
as $u\to\infty$. In this case the value $w_1(x,0)$ (cf.\ \eq{e:wlax0} with $\lambda=1$)
is known as the {\em $\vfi$-variation\/} of the function $x:[a,b]\to\RB$
(in the sense of F.\,Riesz, Yu.\,T.\,Medvedev and W.\,Orlicz), the function $x$ with
$w_1(x,0)<\infty$ is said to be {\em of bounded $\vfi$-variation\/} on $[a,b]$, and
we have:
  \begin{equation} \label{e:x-y}
w_\lambda(x,y)=w_\lambda(x-y,0)=w_1\Biggl(\frac{x-y}\lambda,0\Biggr),
\quad\lambda>0,\!\!\quad x,y\in\XB=\RB^{[a,b]}.
  \end{equation}

Denote by $\mbox{AC}[a,b]$ the space of all absolutely continuous real valued
functions on $[a,b]$ and by $\mbox{L}^1[a,b]$ the space of all (equivalence
classes of) Lebesgue summable functions on~$[a,b]$.

\smallbreak
The following criterion is known for functions $x:[a,b]\to\RB$ to be in the space
$\GV\vfi[a,b]=\XB_w^*$ (for more details see \cite{gens}, \cite[Sections 3, 4]{Sovae},
\cite{Cy-Ma}, \cite[Section~2.4]{Maligranda}, \cite{Medvedev}):
$x\in\GV\vfi[a,b]$ iff $w_\lambda(x,0)=w_1(x/\lambda,0)<\infty$ for some $\lambda=\lambda(x)>0$ (i.e., $x/\lambda$ is of bounded $\vfi$-variation on $[a,b]$)
iff $x\in\mbox{AC}[a,b]$ and its derivative $x'\in\mbox{L}^1[a,b]$
(defined almost everywhere on $[a,b]$) satisfies the condition:
  \begin{equation} \label{e:intab}
w_\lambda(x,x_0)=w_\lambda(x,0)=\int_a^b\vfi\Biggl(\frac{|x'(t)|}\lambda\Biggr)dt
<\infty,\quad\,x_0\in\RB.
  \end{equation}

Given $x_0\in\RB$, we set $X=\{x:[a,b]\to\RB\mid x(a)=x_0\}$, and so (cf.\
\eq{e:GVX}),
  \begin{equation} \label{e:Xax}
X_w^*=X_w^*(x_0)=\{x\in\GV\vfi[a,b]:x(a)=x_0\}.
  \end{equation}
Thus, the modular $w$ is strict and convex on $X$ and the modular space
\eq{e:Xax} is modular complete. Note that $X_w^*$ is {\em not\/} a linear
subspace of $\GV\vfi[a,b]$, which is a normed Banach algebra
(cf.\ \cite[Theorem~3.6]{posit}).

\smallbreak
{\bf\ref{s:exas}.5.} Here we present an example when the metric and modular
convergences coincide. This example is a modification of Example~3.5(c) from~\cite{Sovae}.
We set $[a,b]=[0,1]$, $M=\RB$ and $\vfi(u)=e^u-1$ for $u\in\RB^+$.
Clearly, $\vfi$ satisfies the Orlicz condition, but does not
satisfy the $\Delta_2$-condition at infinity.

\smallbreak
Given a number $\alpha>0$, we define a function $x_\alpha:[0,1]\to\RB$ by
  $$x_\alpha(t)=\alpha t(1-\log t)\quad\mbox{if}\quad 0<t\le1\quad
     \mbox{and}\quad x_\alpha(0)=0.$$
Since $x_\alpha'(t)=-\alpha\log t$ for $0<t\le1$, by \eq{e:intab}, for any
number $\lambda>0$ we find
  $$w_\lambda(x_\alpha,0)=
     \int_0^1\vfi\Biggl(\frac{|x_\alpha'(t)|}\lambda\Biggr)dt=
     \int_0^1\frac{dt}{t^{\alpha/\lambda}}-1=
     \left\{\begin{array}{ccr}
     \infty                                     &\mbox{if}& 0<\lambda\le\alpha,\\[4pt]
     \D\frac{\alpha}{\lambda-\alpha} &\mbox{if}& \lambda>\alpha.
             \end{array}
     \right.$$
It follows that the modular $w$ can take infinite values (although it is strict) and that
$x_\alpha\in X_w^*=X_w^*(0)$ for all $\alpha>0$. Also, we have:
  $$d_w^*(x_\alpha,0)=\inf\{\lambda>0:w_\lambda(x_\alpha,0)\le1\}=2\alpha.$$
Thus, if we set $\alpha=\alpha(n)=1/n$ and $x_n=x_{\alpha(n)}$ for $n\in\NB$,
then we find that $d_w^*(x_n,0)\to0$ as $n\to\infty$ and $w_\lambda(x_n,0)\to0$
as $n\to\infty$ for all $\lambda>0$, and, in accordance with Theorem~\ref{t:dil=w},
these two convergences are equivalent.

\smallbreak
{\bf\ref{s:exas}.6.} Here we expose an example when the modular convergence
is weaker than the metric convergence. Let $[a,b]$, $M$ and $\vfi$ be as in
Example~\ref{s:exas}.5.

\smallbreak
Given $0\le\beta\le1$, we define a function $x_\beta:[0,1]\to\RB$ as follows:
  $$x_\beta(t)=t-(t+\beta)\log(t+\beta)+\beta\log\beta\quad\mbox{if}\quad
     \beta>0\quad\mbox{and}\quad 0\le t\le1$$
and
  $$x_0(t)=t-t\log t\quad\mbox{if}\quad 0<t\le1\quad\mbox{and}\quad x_0(0)=0.$$
Since $x_\beta'(t)=-\log(t+\beta)$ for $\beta>0$ and $t\in[0,1]$, we have:
  $$|x_\beta'(t)|\!=\!-\log(t+\beta)\,\,\mbox{if}\,\,0\le t\le 1-\beta,\,\,\mbox{and}\,\,
     |x_\beta'(t)|=\log(t+\beta)\,\,\mbox{if}\,\,1-\beta\le t\le1,$$
and so, by virtue of \eq{e:intab}, given $\lambda>0$, we find
  $$w_\lambda(x_\beta,0)=\int_0^1\vfi(|x_\beta'(t)|/\lambda)\,dt=I_1+I_2-1,
     \qquad\beta>0,$$
where
  $$I_1=\int_0^{1-\beta}\frac{dt}{(t+\beta)^{1/\lambda}}=
     \left\{\begin{array}{ccr}
     \D\frac{\lambda}{\lambda-1}\biggl(1-\beta^{(\lambda-1)/\lambda}\biggr)
                   &\mbox{if}& 0<\lambda\ne1,\\[4pt]
     -\log\beta &\mbox{if}& \lambda=1,
             \end{array}
     \right.$$
and
  $$I_2=\int_{1-\beta}^1(t+\beta)^{1/\lambda}\,dt=
     \frac{\lambda}{\lambda+1}\biggl((1+\beta)^{(\lambda+1)/\lambda}-1\biggr)
     \quad\mbox{for all}\quad\lambda>0.$$

Also, $w_\lambda(x_0,0)=\infty$ if $0<\lambda\le1$, and
$w_\lambda(x_0,0)=1/(\lambda-1)$ if $\lambda>1$
(cf.\ Example~\ref{s:exas}.5 with $\alpha=1$). Thus,
$x_\beta\in X_w^*=X_w^*(0)$ for all $0\le\beta\le1$.

\smallbreak
Clearly, $x_\beta$ converges pointwise on $[0,1]$ to $x_0$ as $\beta\to+0$
(actually, the first inequality in the proof of \cite[Lemma~4.1(a)]{Sovae} shows
that the convergence is uniform on $[0,1]$).

\smallbreak
Now we calculate the values $w_\lambda(x_\beta,x_0)$ for $\lambda>0$ and
$d_w^*(x_\beta,x_0)$ and invesigate their convergence to zero as $\beta\to+0$. Since
  $$(x_\beta-x_0)'(t)=-\log(t+\beta)+\log t\quad\mbox{for}\quad 0<t\le1,$$
we have:
  $$\frac{|(x_\beta-x_0)'(t)|}{\lambda}=\frac{\log(t+\beta)-\log t}{\lambda}=
     \log\biggl(1+\frac\beta t\biggr)^{1/\lambda},$$
and so, by virtue of \eq{e:x-y} and \eq{e:intab},
  $$w_\lambda(x_\beta,x_0)=
     \int_0^1\vfi\Biggl(\frac{|(x_\beta-x_0)'(t)|}\lambda\Biggr)dt=
     -1+\int_0^1\biggl(1+\frac\beta t\biggr)^{1/\lambda}dt.$$

If $0<\lambda\le1$, we have
  $$\biggl(1+\frac\beta t\biggr)^{1/\lambda}\ge 1+\frac\beta t\quad\,\,\mbox{and}
     \quad\,\,\int_0^1\biggl(1+\frac\beta t\biggr)dt=\infty,$$
and so, $w_\lambda(x_\beta,x_0)=\infty$ for all $0<\beta\le1$ and $0<\lambda\le1$.

\smallbreak
Now suppose that $\lambda>1$. Then
  $$w_\lambda(x_\beta,x_0)=-1+\int_0^\beta\biggl(1+\frac\beta t\biggr)^{1/\lambda}dt
     +\int_\beta^1\biggl(1+\frac\beta t\biggr)^{1/\lambda}dt\equiv-1+I\!I_1+I\!I_2,$$
where
  \begin{eqnarray}
I\!I_1\!\!&\le\!\!&\D\int_0^\beta\biggl(\frac{2\beta}{t}\biggr)^{1/\lambda}dt=
  (2\beta)^{1/\lambda}\!\int_0^\beta t^{-1/\lambda}dt=
  (2\beta)^{1/\lambda}\!\cdot\!\frac\lambda{\lambda-1}\!\cdot\!%
  \beta^{1-(1/\lambda)}=\nonumber\\[2pt]
\!\!&=\!\!&2^{1/\lambda}\!\cdot\!\D\frac{\lambda\beta}{\lambda-1}\to0\quad
  \mbox{as}\quad\beta\to+0\nonumber
  \end{eqnarray}
and
  $$I\!I_2\le\int_\beta^1\biggl(1+\frac\beta t\biggr)dt=
     (1-\beta)-\beta\log\beta\to1\quad\mbox{as}\quad\beta\to+0.$$
It follows that $w_\lambda(x_\beta,x_0)\to0$ as $\beta\to+0$ for all $\lambda>1$.

\smallbreak
On the other hand, since $w_\lambda(x_\beta,x_0)=\infty$ for all $0<\beta\le1$
and $0<\lambda\le1$ (as noticed above), we get
  $d_w^*(x_\beta,x_0)=\inf\{\lambda>0:w_\lambda(x_\beta,x_0)\le1\}\ge1$,
and so, $d_w^*(x_\beta,x_0)$ cannot converge to zero as $\beta\to+0$.

\smallbreak
Thus, if we set $\beta=\beta(n)=1/n$ and $x_n=x_{\beta(n)}$ for $n\in\NB$, then
we find $d_w^*(x_n,x_0)\not\to0$ as $n\to\infty$, whereas $w_\lambda(x_n,x_0)\to0$
as $n\to\infty$ only for $\lambda>1$.

\section{A fixed point theorem for modular contractions} \label{s:fp}

Since convex modulars play the central role in this section, we concentrate mainly
on them. We begin with a characterization of $d_w^*$-Lipschitz maps on the
modular space $X_w^*$ in terms of their generating convex modulars~$w$.

\begin{thm} \label{t:Lip}
Let $w$ be a convex modular on $X$ and $k>0$ be a constant. Given a map
$T:X_w^*\to X_w^*$ and $x,y\in X_w^*$, the Lipschitz condition
$d_w^*(Tx,Ty)\le k\,d_w^*(x,y)$ is equivalent to the following\/{\rm:}
$w_{k\lambda+0}(Tx,Ty)\le1$ for all $\lambda>0$ such that $w_\lambda(x,y)\le1$.
\end {thm}

\begin{pf}
First, note that, given $c>0$, the function, defined by
$\ov w_\lambda(x,y)=w_{c\lambda}(x,y)$, $\lambda>0$, $x,y\in X$, is also a convex
modular on $X$ and $d_{\ov w}^*=\frac1c d_w^*$:
  \begin{eqnarray}
d_{\ov w}^*(x,y)\!\!&=\!\!&\inf\{\lambda\!>\!0:w_{c\lambda}(x,y)\!\le\!1\}\!=\!
     \inf\{\mu/c\!>\!0:w_\mu(x,y)\!\le\!1\}=\nonumber\\
\!\!&=\!\!&\frac1c\,d_w^*(x,y)\quad\,\,\mbox{for \,all \,$x,y\in X_{\ov w}^*=X_w^*$}.
  \label{e:1c}
  \end{eqnarray}

\smallbreak
{\em Necessity}. We may suppose that $x\ne y$. For any $c>k$, by the
assumption, we find $d_w^*(Tx,Ty)\le k\,d_w^*(x,y)<c\,d_w^*(x,y)$, whence
$d_w^*(Tx,Ty)/c<d_w^*(x,y)$. It follows that if $\lambda>0$ is such that
$w_\lambda(x,y)\le1$, then, by \eq{e:dw*}, $d_w^*(x,y)\le\lambda$ implying,
in view of \eq{e:1c},
  $$\lambda>\frac1c\,d_w^*(Tx,Ty)=\inf\{\mu>0:w_{c\mu}(Tx,Ty)\le1\},$$
and so, $w_{c\lambda}(Tx,Ty)\le1$. Passing to the limit as $c\to k+0$, we arrive
at the desired inequality $w_{k\lambda+0}(Tx,Ty)\le1$.

\smallbreak
{\em Sufficiency}. By the assumption, the set $\{\lambda>0:w_\lambda(x,y)\le1\}$
is contained in the set $\{\lambda>0:w_{k\lambda}^+(Tx,Ty)=%
w_{k\lambda+0}(Tx,Ty)\le1\}$, and so, taking the infima, by virtue of
\eq{e:dw*}, \eq{e:1c} and the equality $d_{w^+}^*=d_w^*$, we get
  $$d_w^*(x,y)\ge\frac1k\,d_{w^+}^*(Tx,Ty)=\frac1k\,d_w^*(Tx,Ty),$$
which implies that $T$ satisfies the Lipschitz condition with constant~$k$.
\qed\end{pf}

Theorem~\ref{t:Lip} can be reformulated as follows. Since
(cf.\ \cite[Theorem~3.8(a)]{NA-I} and \eq{e:dw*}), for $\lambda^*=d_w^*(x,y)$,
  $$(\lambda^*,\infty)\subset\{\lambda\!>\!0:w_\lambda(x,y)\!<\!1\}\subset
     \{\lambda\!>\!0:w_\lambda(x,y)\!\le\!1\}\subset[\lambda^*,\infty),$$
we have: $d_w^*(Tx,Ty)\!\le\! k\,d_w^*(x,y)$ iff $w_{k\lambda}(Tx,Ty)\le1$
for all $\lambda\!>\!\lambda^*\!=\!d_w^*(x,y)$.

\smallbreak
For a metric space $(X,d)$ and the modular $w$ from \eq{e:dxy} on it,
Theorem~\ref{t:Lip} gives the usual Lipschitz condition:
$d(Tx,Ty)/(k\lambda)=w_{k\lambda}(Tx,Ty)\le1$ for all $\lambda>0$ such
that $d(x,y)/\lambda=w_\lambda(x,y)\le1$, i.e., $d(Tx,Ty)\le k\lambda$ for
all $\lambda\ge d(x,y)$, and so, $d(Tx,Ty)\le kd(x,y)$.

\smallbreak
As a corollary of Theorem~\ref{t:Lip}, we find that
  \begin{equation} \label{e:corLip}
\mbox{if $w_{k\lambda}(Tx,Ty)\!\le\!w_\lambda(x,y)$ for all $\lambda\!>\!0$, then
$d_w^*(Tx,Ty)\!\le\! k\,d_w^*(x,y)$;}
  \end{equation}
in fact, it suffices to note only that if $\lambda>0$ is such that $w_\lambda(x,y)\le1$,
then, by \eq{e:lapm}, $w_{k\lambda+0}(Tx,Ty)\le w_{k\lambda}(Tx,Ty)\le%
w_\lambda(x,y)\le1$, and apply Theorem~\ref{t:Lip}.

\smallbreak
Now we briefly comment on $d_w$-Lipschitz maps on $X_w^*$, where $w$ is a
general modular on~$X$ and $d_w$ is the metric from \eq{e:dw}.
Note that, given $c>0$, the function $\ov w_\lambda(x,y)=\frac1c\,w_{c\lambda}(x,y)$
is also a modular on $X$ and $d_{\ov w}=\frac1c\,d_w$ on $X_{\ov w}^*=X_w^*$.
Following the lines of the proof of Theorem~\ref{t:Lip}, we get

\begin{thm} \label{t:Lip2}
If $w$ is a modular on $X$ and $k>0$, given $T:X_w^*\to X_w^*$ and $x,y\in X_w^*$,
we have\/{\rm} $d_w(Tx,Ty)\le k\,d_w(x,y)$ iff\/ $w_{k\lambda+0}(Tx,Ty)\le k\lambda$
for all $\lambda>0$ such that $w_\lambda(x,y)\le\lambda$.
\end{thm}

The following assertion is a corollary of Theorem~\ref{t:Lip2}:
  $$\mbox{if $w_{k\lambda}(Tx,Ty)\!\le\!k\,w_\lambda(x,y)$ for all $\lambda\!>\!0$,
     then $d_w(Tx,Ty)\!\le\! k\,d_w(x,y)$.}$$

\begin{defin} \label{d:modcontr}
Given a (convex) modular $w$ on $X$, a map $T:X_w^*\to X_w^*$ is said to be
{\em modular contractive\/} (or a {\em $w$-contraction\/}) provided there exist 
numbers $0<k<1$ and $\lambda_0>0$, possibly depending on~$k$, such that
  \begin{equation} \label{e:mc}
\mbox{$w_{k\lambda}(Tx,Ty)\le w_\lambda(x,y)$ \,\,for \,all \,$0<\lambda\le\lambda_0$
\,and \,$x,y\in X_w^*$.}
  \end{equation}
\end{defin}

A few remarks are in order. First, by virtue of \eq{e:dxy}, for a metric space
$(X,d)$ condition \eq{e:mc} is equivalent to the usual one:
$d(Tx,Ty)\le kd(x,y)$. Second, condition \eq{e:mc} is a {\em local\/} one
with respect to $\lambda$ as compared to the assumption on the left in \eq{e:corLip},
and the principal inequality in it may be of the form $\infty\le\infty$.
Third, if, in addition, $w$ is {\em strict\/} and if we set $\infty/\infty=1$, then \eq{e:mc}
is a consequence of the following: there exists a number $0\!<\!h\!<\!1$ such that
  \begin{equation} \label{e:mc2}
\limsup_{\lambda\to+0}\,\Biggl(\sup_{x\ne y}\frac{w_{h\lambda}(Tx,Ty)}%
{w_\lambda(x,y)}\Biggr)\le1,
  \end{equation}
where the supremum is taken over all $x,y\in X_w^*$ such that $x\ne y$.
In order to see this, we first note that the left hand side in \eq{e:mc2} is well
defined in the sense that, by virtue of (\is) from definition~\ref{d:md},
$w_\lambda(x,y)\ne0$ for all $\lambda>0$ and $x\ne y$. Choose any
$k$ such that $h<k<1$. It follows from \eq{e:mc2} that
  $$\lim_{\mu\to+0}\,\sup_{\lambda\in(0,\mu]}\,\Biggl(
     \sup_{x\ne y}\frac{w_{h\lambda}(Tx,Ty)}{w_\lambda(x,y)}\Biggr)\le1<\frac kh,$$
and so, there exists a $\mu_0=\mu_0(k)>0$ such that
  $$\sup_{x\ne y}\frac{w_{h\lambda}(Tx,Ty)}{w_\lambda(x,y)}<\frac kh
     \quad\,\,\mbox{for \,all}\quad\,\, 0<\lambda\le\mu_0,$$
whence
  $$w_{h\lambda}(Tx,Ty)\le\frac kh\,w_\lambda(x,y),\quad\,\, 0<\lambda\le\mu_0,
     \quad\,\, x,y\in X_w^*.$$
Taking into account inequalities \eq{e:mpc} and $(h/k)\lambda<\lambda$, we get
  $$w_\lambda(x,y)\le\frac{(h/k)\lambda}{\lambda}\,w_{(h/k)\lambda}(x,y)=
     \frac hk\,w_{(h/k)\lambda}(x,y),$$
which together with the previous inequality gives:
  $$w_{h\lambda}(Tx,Ty)\le w_{(h/k)\lambda}(x,y)\quad\mbox{for all}
     \quad 0<\lambda\le\mu_0\quad
    \mbox{and}\quad x,y\in X_w^*.$$
Setting $\lambda'\!=\!(h/k)\lambda$ and $\lambda_0\!=\!(h/k)\mu_0$ and noting that
$0\!<\!\lambda'\!\le\!\lambda_0$ and $h\lambda\!=\!k\lambda'$, the last inequality implies
$w_{k\lambda'}(Tx,Ty)\le w_{\lambda'}(x,y)$ for all $0<\lambda'\le\lambda_0$ and
$x,y\in X_w^*$, which is exactly \eq{e:mc}.

\smallbreak
The main result of this paper is the following fixed point theorem for modular
contractions in modular metric spaces~$X_w^*$.

\begin{thm} \label{t:main}
Let $w$ be a strict convex modular on $X$ such that the modular space $X_w^*$
is $w$-complete, and $T:X_w^*\to X_w^*$ be a $w$-contractive map such that
  \begin{equation} \label{e:reach}
\mbox{for each $\lambda\!>\!0$ there exists an $x\!=\!x(\lambda)\!\in\! X_w^*$
such that  $w_\lambda(x,Tx)\!<\!\infty$.}
  \end{equation} 
Then $T$ has a fixed point, i.e., $Tx_*=x_*$ for some $x_*\in X_w^*$.
If, in addition, the modular $w$ assumes only finite values on $X_w^*$, then condition\/
\eq{e:reach} is redundant, the fixed point $x_*$ of\/ $T$ is unique and for each
$\ov x\in X_w^*$ the sequence of iterates $\{T^n\ov x\}$ is modular
convergent to~$x_*$. 
\end{thm}

\begin{pf}
Since $w$ is convex, the following inequality follows by induction from
condition (iv) of definition~\ref{d:md}:
  \begin{equation} \label{e:cla}
(\lambda_1+\lambda_2+\cdots+\lambda_N)w_{\lambda_1+\lambda_2+\cdots+\lambda_N}
(x_1,x_{N+1})\le\sum_{i=1}^N\lambda_i w_{\lambda_i}(x_i,x_{i+1}),
  \end{equation}
where $N\in\NB$, $\lambda_1,\lambda_2,\dots,\lambda_N\in(0,\infty)$
and $x_1,x_2,\dots,x_{N+1}\in X$. In the proof below we will need a variant
of this inequality. Let $n,m\in\NB$, $n>m$,
$\lambda_m,\lambda_{m+1},\dots,\lambda_{n-1}\in(0,\infty)$ and
$x_m,x_{m+1},\dots,x_n\in X$. Setting $N=n-m$,
$\lambda_j'=\lambda_{j+m-1}$ for $j=1,2,\dots,N$, $x_j'=x_{j+m-1}$ for
$j=1,2,\dots,N+1$ and applying \eq{e:cla} to the primed lambda's and $x$'s, we get:
  \begin{equation} \label{e:varcla}
(\lambda_m\!+\!\lambda_{m+1}\!+\!\cdots\!+\!\lambda_{n-1})
w_{\lambda_m+\lambda_{m+1}+\cdots+\lambda_{n-1}}(x_m,x_n)\le
\sum_{i=m}^{n-1}\lambda_iw_{\lambda_i}(x_i,x_{i+1}).
  \end{equation}

By the $w$-contractivity of $T$, there exist two numbers $0<k<1$ and
$\lambda_0=\lambda_0(k)>0$ such that condition \eq{e:mc} holds. Setting
$\lambda_1=(1-k)\lambda_0$, the  assumption \eq{e:reach} implies the existence
of an element $\ov x=\ov x(\lambda_1)\in X_w^*$ such that
$C=w_{\lambda_1}(\ov x,T\ov x)$ is finite. We set $x_1=T\ov x$ and
$x_n=Tx_{n-1}$ for all integer $n\ge2$, and so, $\{x_n\}\subset X_w^*$ and
$x_n=T^n\ov x$, where $T^n$ designates the $n$-th iterate of~$T$.
We are going to show that the sequence $\{x_n\}$ is $w$-Cauchy.
Since $k^i\lambda_1<\lambda_1<\lambda_0$ for all $i\in\NB$, inequality \eq{e:mc} yields:
  $$w_{k^i\lambda_1}(x_i,x_{i+1})=w_{k(k^{i-1}\lambda_1)}(Tx_{i-1},Tx_i)\le
     w_{k^{i-1}\lambda_1}(x_{i-1},x_i),$$
and it follows by induction that
  \begin{equation} \label{e:wki}
w_{k^i\lambda_1}(x_i,x_{i+1})\le w_{\lambda_1}(\ov x,x_1)=C\quad
\mbox{for \,all}\quad i\in\NB.
  \end{equation}
Let integers $n$ and $m$ be such that $n>m$. We set
  $$\lambda=\lambda(n,m)=k^m\lambda_1+k^{m+1}\lambda_1+\cdots+
     k^{n-1}\lambda_1=k^m\frac{1-k^{n-m}}{1-k}\,\lambda_1.$$
By virtue of \eq{e:varcla} with $\lambda_i=k^i\lambda_1$ and \eq{e:wki}, we find
  $$w_\lambda(x_m,x_n)\le\sum_{i=m}^{n-1}\frac{k^i\lambda_1}\lambda\,%
     w_{k^i\lambda_1}(x_i,x_{i+1})\le\frac1\lambda\,%
     \Biggl(\sum_{i=m}^{n-1}k^i\lambda_1\Biggr)C=C,\quad n>m.$$
Taking into account that
  $$\lambda_0=\frac{\lambda_1}{1-k}>k^m\frac{1-k^{n-m}}{1-k}\,\lambda_1
     =\lambda(n,m)=\lambda\quad\mbox{for \,all}\quad n>m,$$
and applying \eq{e:mpc}, we get:
  $$w_{\lambda_0}(x_m,x_n)\le\frac{\lambda}{\lambda_0}\,w_\lambda(x_m,x_n)
     \le k^m\frac{1\!-\!k^{n-m}}{1-k}\!\cdot\!\frac{\lambda_1}{\lambda_0}\,C
     \le k^mC\to0\,\,\mbox{as}\,\,m\to\infty.$$
Thus, the sequence $\{x_n\}$ is modular Cauchy, and so, by the $w$-completeness
of $X_w^*$, there exists an $x_*\in X_w^*$ such that
  $$w_{\lambda_0}(x_n,x_*)\to0\quad\,\,\mbox{as}\quad\,\,n\to\infty.$$
Since $w$ is strict, by Theorem~\ref{t:corr}(b), the modular limit $x_*$ of the
sequence $\{x_n\}$ is determined uniquely.

\smallbreak
Let us show that $x_*$ is a fixed point of $T$, i.e., $Tx_*=x_*$. In fact, by
property (iii) of definition~\ref{d:md} and \eq{e:mc}, we have (note that
$Tx_n=x_{n+1}$):
  \begin{eqnarray}
w_{(k+1)\lambda_0}(Tx_*,x_*)\!\!&\le\!\!&w_{k\lambda_0}(Tx_*,Tx_n)+
  w_{\lambda_0}(x_*,x_{n+1})\le\nonumber\\[4pt]
\!\!&\le\!\!&w_{\lambda_0}(x_*,x_n)+w_{\lambda_0}(x_*,x_{n+1})
  \to0\quad\mbox{as}\quad n\to\infty,\nonumber
  \end{eqnarray}
and so, $w_{(k+1)\lambda_0}(Tx_*,x_*)=0$. By the strictness of $w$, $Tx_*=x_*$.

\smallbreak
Finally, assuming $w$ to be finite valued on $X_w^*$, we show that the fixed point
of $T$ is unique. Suppose $x_*,y_*\in X_w^*$ are such that $Tx_*=x_*$ and
$Ty_*=y_*$. Then the convexity of $w$ and inequalities $k\lambda_0<\lambda_0$
and \eq{e:mc}~imply
  $$w_{\lambda_0}(x_*,y_*)\le\frac{k\lambda_0}{\lambda_0}\,w_{k\lambda_0}%
     (x_*,y_*)=kw_{k\lambda_0}(Tx_*,Ty_*)\le kw_{\lambda_0}(x_*,y_*),$$
and since $w_{\lambda_0}(x_*,y_*)$ is finite, $(1-k)w_{\lambda_0}(x_*,y_*)\le0$.
Thus, $w_{\lambda_0}(x_*,y_*)=0$, and by the strictness of $w$, we get $x_*=y_*$.
The last assertion is clear.
\qed\end{pf}

It is to be noted that assumption \eq{e:reach} in Theorem~\ref{t:main} is (probably)
too strong, and what we actually need for the iterative procedure to work in the
proof of Theorem~\ref{t:main} is only the existence of an $\ov x\in X_w^*$ such
that $w_{(1-k)\lambda_0}(\ov x,T\ov x)<\infty$, where $\lambda_0$ is the
constant from~\eq{e:mc}.

\smallbreak
A standard corollary of Theorem~\ref{t:main} is as follows: if $w$ is finite valued on
$X_w^*$ and an $n$-th iterate $T^n$ of $T:X_w^*\to X_w^*$ satisfies the
assumptions of Theorem~\ref{t:main}, then $T$ has a unique fixed point.
In fact, by Theorem~\ref{t:main} applied to $T^n$, $T^nx_*=x_*$ for some
$x_*\in X_w^*$. Since $T^n(Tx_*)=T(T^nx_*)=Tx_*$, the point $Tx_*$ is also a
fixed point of $T^n$, and so, the uniqueness of a fixed point of $T^n$ implies
$Tx_*=x_*$. We infer that $x_*$ is a unique fixed point of $T$: if $y_*\in X_w^*$
and $Ty_*=y_*$, then $T^ny_*=T^{n-1}(Ty_*)=T^{n-1}y_*=\dots=y_*$,
i.e., $y_*$ is yet another fixed point of $T^n$, and again the uniqueness
of a fixed point of $T^n$ yields $y_*=x_*$.

\smallbreak
Another corollary of Theorem~\ref{t:main} concerns general (nonconvex) modulars
$w$ on~$X$ (cf.\ Theorem~\ref{t:maineconv}).
Taking into account Theorem~\ref{t:Lip2} and its corollary, we have

\begin{defin} \label{d:mcnon}
Given a modular $w$ on $X$, a map $T:X_w^*\to X_w^*$ is said to be {\em strongly 
modular contractive\/} (or a {\em strong $w$-contraction\/}) if there exist 
numbers $0<k<1$ and $\lambda_0=\lambda_0(k)>0$ such that
  \begin{equation} \label{e:mcn}
\mbox{$w_{k\lambda}(Tx,Ty)\le kw_\lambda(x,y)$ \,\,for \,all \,$0<\lambda\le\lambda_0$
\,and \,$x,y\in X_w^*$.}
  \end{equation}
\end{defin}

Clearly, condition \eq{e:mcn} implies condition \eq{e:mc}.

\begin{thm} \label{t:maineconv}
Let $w$ be a strict modular on $X$ such that $X_w^*$ is $w$-complete,
and $T:X_w^*\to X_w^*$ be a strongly $w$-contractive map such that
condition\/ \eq{e:reach} holds. Then $T$ admits a fixed point.
If, in addition, $w$ is finite valued on $X_w^*$, then\/
\eq{e:reach} is redundant, the fixed point $x_*$ of\/ $T$ is unique and for each
$\ov x\in X_w^*$ the sequence of iterates $\{T^n\ov x\}$ is modular
convergent to~$x_*$. 
\end{thm}

\begin{pf}
We set $v_\lambda(x,y)=w_\lambda(x,y)/\lambda$ for all $\lambda>0$ and $x,y\in X$.
It was observed in Section~\ref{s:mm} that $v$ is a convex modular on~$X$.
It is also clear that $v$ is strict and the modular space $X_v^*=X_w^*$ is
$v$-complete. Moreover, condition \eq{e:mcn} for $w$ implies condition \eq{e:mc}
for~$v$, and \eq{e:reach} is satisfied with $w$ replaced by~$v$. By Theorem~\ref{t:main},
applied to~$X$ and $v$, there exists an $x_*\in X_v^*=X_w^*$ such that
$Tx_*=x_*$. The remaining assertions are obvious.
\qed\end{pf}

\section{An application of the fixed point theorem} \label{s:appl}

In this section we present a rather standard application of Theorem~\ref{t:main}
to the Carath{\'e}odory-type ordinary differential equations. The key interest will
be in obtaining the inequality~\eq{e:mc}.

\smallbreak
Given a convex $\vfi$-function $\vfi$ on $\RB^+$ satisfying the Orlicz condition
at infinity, we denote by $\mbox{L}^\vfi[a,b]$ the Orlicz space of real valued functions
on $[a,b]$ (cf.\ \cite[Chapter~II]{Musielak}),
i.e., a function $z:[a,b]\to\RB$ (or an almost everywhere finite valued function $z$
on $[a,b]$) belongs to $\mbox{L}^\vfi[a,b]$ provided $z$ is measurable and
$\rho(z/\lambda)<\infty$ for some number $\lambda=\lambda(z)>0$, where
$\rho(z)=\int_a^b\vfi(|z(t)|)dt$ is the classical Orlicz modular.

\smallbreak
Suppose $f:[a,b]\times\RB\to\RB$ is a (Carath{\'e}odory-type) function, which
satisfies the following two conditions:

\medbreak
(C.1) for each $x\in\RB$ the function $f(\cdot,x)=[t\mapsto f(t,x)]$ is
measurable on $[a,b]$ and there exists a point $y_0\in\RB$ such that
$f(\cdot,y_0)\in\mbox{L}^\vfi[a,b]$;

\smallbreak
(C.2) there exists a constant $L>0$ such that $|f(t,x)-f(t,y)|\le L|x-y|$
for almost all $t\in[a,b]$ and all $x,y\in\RB$.

\medbreak
Given $x_0\in\RB$, we let $X_w^*$ be the modular space \eq{e:Xax} generated
by the modular $w$ from \eq{e:wlaxy} under the assumptions from
Section~\ref{s:exas}.4.

\smallbreak
Consider the following integral operator
  \begin{equation} \label{e:intop}
(Tx)(t)=x_0+\int_a^t\! f(s,x(s))\,ds,\qquad x\in X_w^*,\quad t\in[a,b].
  \end{equation}

\begin{thm} \label{t:meniq}
Under the assumptions\/ {\rm(C.1)} and\/  {\rm(C.2)} the operator $T$ maps
$X_w^*$ into itself, and the following inequality holds  in $[0,\infty]${\rm:}
  \begin{equation} \label{e:Lba}
w_{L(b-a)\lambda}(Tx,Ty)\le w_\lambda(x,y)\quad\mbox{for \,all \,$\lambda>0$
 \,and \,$x,y\in X_w^*$.}
  \end{equation}
\end{thm}

\begin{pf}
We will apply the {\em Jensen integral inequality\/} with the convex $\vfi$-fun\-c\-tion
$\vfi$ (e.g., \cite[X.5.6]{Natanson}) several times:
  \begin{equation} \label{e:Jensen}
\vfi\Biggl(\frac1{b-a}\int_a^b|x(t)|dt\biggr)\le\frac1{b-a}\int_a^b
\vfi\bigl(|x(t)|\bigr)dt,\quad\,\,x\in\mbox{L}^1[a,b],
  \end{equation}
where the intergral in the right hand side is well defined in the sense that it takes
values in~$[0,\infty]$.

\smallbreak
1. First, we show that $T$ is well defined on $X_w^*$. Let $x\in X_w^*$, i.e.,
$x\in\GV\vfi[a,b]$ and $x(a)=x_0$. Since (cf.\ Section~\ref{s:exas}.4)
$x\in\mbox{AC}[a,b]$, by virtue of (C.1) and (C.2), the composed function
$t\mapsto f(t,x(t))$ is measurable on~$[a,b]$. Let us prove that this function
belongs to $\mbox{L}^1[a,b]$. By Lebesgue's Theorem,
$x(t)=x_0+\int_a^t x'(s)ds$ for all $t\in[a,b]$, and so, (C.2) yields
  \begin{eqnarray}
|f(t,x(t))|\!\!&\le\!\!&|f(t,x(t))-f(t,y_0)|+|f(t,y_0)|\le\nonumber\\[6pt]
\!\!&\le\!\!&L|x(t)-y_0|+|f(t,y_0)|\le\nonumber\\
\!\!&\le\!\!&\D L\int_a^b|x'(s)|ds+L|x_0-y_0|+|f(t,y_0)|\label{e:6.3}
  \end{eqnarray}
for almost all $t\in[a,b]$. Since $x\in X_w^*$, and so, $x\in\GV\vfi[a,b]$,
there exists a constant $\lambda_1=\lambda_1(x)>0$ such that (cf.\ \eq{e:intab})
  $$C_1\equiv w_{\lambda_1}(x,x_0)=
     \int_a^b\vfi\Biggl(\frac{|x'(s)|}{\lambda_1}\Biggr)ds<\infty,$$
and since, by (C.1), $f(\cdot,y_0)\!\in\!\mbox{L}^\vfi[a,b]$, there exists a
constant \mbox{$\lambda_2\!=\!\lambda_2(f(\cdot,y_0))\!>\!0$} such that
  $$C_2\equiv\rho\bigl(f(\cdot,y_0)/\lambda_2\bigr)=
    \int_a^b\vfi\Biggl(\frac{|f(t,y_0)|}{\lambda_2}\Biggr)dt<\infty.$$
Setting $\lambda_0=L\lambda_1(b-a)+1+\lambda_2$ and noting that
  $$\frac{L\lambda_1(b-a)}{\lambda_0}+\frac{1}{\lambda_0}+
     \frac{\lambda_2}{\lambda_0}=1,$$
by the convexity of $\vfi$, we find (see \eq{e:6.3})
  \begin{eqnarray}
&&\qquad\quad\vfi\Biggl(\frac1{\lambda_0}\,\Biggl[L\int_a^b|x'(s)|ds+
  L|x_0-y_0|+|f(t,y_0)|\Biggr]\Biggr)\le\nonumber\\[4pt]
&\le\!\!&\frac{L\lambda_1(b\!-\!a)}{\lambda_0}\,\vfi\Biggl(
  \frac1{b\!-\!a}\int_a^b\!\frac{|x'(s)|}{\lambda_1}\,ds\Biggr)\!+\!
  \frac1{\lambda_0}\,\vfi\bigl(L|x_0\!-\!y_0|\bigr)\!+\!
  \frac{\lambda_2}{\lambda_0}\,\vfi\Biggl(\frac{|f(\cdot,y_0)|}{\lambda_2}\Biggr),
  \nonumber
  \end{eqnarray}
and so, \eq{e:6.3} and Jensen's integral inequality yield
  \begin{equation} \label{e:abfa}
\int_a^b\!\vfi\Biggl(\frac{|f(t,x(t))|}{\lambda_0}\Biggr)dt\le
\frac{L\lambda_1(b\!-\!a)}{\lambda_0}C_1+
\frac{b\!-\!a}{\lambda_0}\vfi\bigl(L|x_0-y_0|\bigr)+
\frac{\lambda_2}{\lambda_0}C_2\equiv C_0<\infty.
  \end{equation}
Now, it follows from \eq{e:Jensen} that
  $$\vfi\Biggl(\frac1{\lambda_0(b-a)}\int_a^b|f(t,x(t))|dt\Biggr)\le
     \frac1{b-a}\int_a^b\vfi\Biggl(\frac{|f(t,x(t))|}{\lambda_0}\Biggr)dt\le
     \frac{C_0}{b-a}$$
implying
 $$\int_a^b|f(t,x(t))|dt\le\lambda_0(b-a)\vfi^{-1}\Biggl(\frac{C_0}{b-a}\Biggr)<\infty.$$
Thus, $[t\mapsto f(t,x(t))]\in\mbox{L}^1[a,b]$. As a consequence, the operator
$T$ is well defined on $X_w^*$ and, by \eq{e:intop}, $Tx\in\mbox{AC}[a,b]$
for all $x\in X_w^*$, which implies that the almost everywhere derivative
$(Tx)'$ belongs to $\mbox{L}^1[a,b]$ and satisfies
  \begin{equation} \label{e:aetxt}
(Tx)'(t)=f(t,x(t))\quad\mbox{for \,almost \,all}\quad t\in[a,b].
  \end{equation}

2. It is clear from \eq{e:intop} that, given $x\in X_w^*$, $(Tx)(a)=x_0$,
and so, $Tx\in X=\{y:[a,b]\to\RB\mid y(a)=x_0\}$. Now we show that $Tx\in X_w^*$.
In fact, by virtue of \eq{e:intab}, \eq{e:aetxt} and \eq{e:abfa}, we have
  \begin{equation} \label{e:Txx}
w_{\lambda_0}(Tx,x_0)=\int_a^b\vfi\Biggl(\frac{|(Tx)'(t)|}{\lambda_0}\Biggr)dt=
\int_a^b\vfi\Biggl(\frac{|f(t,x(t))|}{\lambda_0}\Biggr)dt\le C_0,
  \end{equation}
and so, $T$ maps $X_w^*$ into itself.

\smallbreak
3. In order to obtain inequality \eq{e:Lba}, let $\lambda>0$ and $x,y\in X_w^*$.
Taking into account \eq{e:x-y}, \eq{e:intab} and \eq{e:aetxt}, we find
  \begin{eqnarray}
w_{L(b-a)\lambda}(Tx,Ty)\!\!&=\!\!&w_{L(b-a)\lambda}(Tx\!-\!Ty,x_0)=
  \int_a^b\vfi\Biggl(\frac{|(Tx\!-\!Ty)'(t)|}{L(b-a)\lambda}\Biggr)dt=\nonumber\\
\!\!&=\!\!&\int_a^b\vfi\Biggl(\frac{|f(t,x(t))-f(t,y(t))|}{L(b-a)\lambda}\Biggr)dt.
  \label{e:wL}
  \end{eqnarray}
Applying (C.2) and Lebesgue's Theorem, we get, for almost all $t\in[a,b]$ (note
that $x(a)=y(a)=x_0$),
  $$|f(t,x(t))-f(t,y(t))|\le L|x(t)-y(t)|\le L\int_a^b|(x-y)'(s)|ds,$$
and so, by \eq{e:Jensen}, the monotonicity of $\vfi$, \eq{e:intab} and \eq{e:x-y},
  \begin{eqnarray}
\vfi\Biggl(\frac{|f(t,x(t))-f(t,y(t))|}{L(b-a)\lambda}\Biggr)\!\!&\le\!\!&
  \vfi\Biggl(\frac1{b-a}\int_a^b\frac{|(x-y)'(s)|}{\lambda}\,ds\Biggr)\le\nonumber\\
\!\!&\le\!\!&\frac1{b-a}\int_a^b\vfi\Biggl(\frac{|(x-y)'(s)|}{\lambda}\Biggr)ds=
  \nonumber\\
\!\!&=\!\!&\frac{1}{b-a}\,w_\lambda(x,y).\nonumber
  \end{eqnarray}
Now, inequality \eq{e:Lba} follows from \eq{e:wL}.
\qed\end{pf}

As a corollary of Theorems \ref{t:main} and \ref{t:meniq}, we have

\begin{thm} \label{t:Cde}
Under the conditions\/ {\rm(C.1)} and\/ {\rm(C.2)}, given $x_0\in\RB$, the
initial value problem
  \begin{equation} \label{e:ivp}
\mbox{$x'(t)=f(t,x(t))$ \,\,for almost all \,\,$t\in[a,b_1]$ \,\,and \,\,$x(a)=x_0$}
  \end{equation}
admits a solution $x\in\GV\vfi[a,b_1]$ with $a<b_1\in\RB$ such that $L(b_1-a)<1$.
\end{thm}

\begin{pf}
We know from Section~\ref{s:exas}.4 that $w$ is a strict convex modular
on the set $X=\{x:[a,b_1]\to\RB\mid x(a)=x_0\}$ and that the modular space
$X_w^*=\GV\vfi[a,b_1]\cap X$ is $w$-complete. By Theorem~\ref{t:meniq},
the operator $T$ from \eq{e:intop} maps $X_w^*$ into itself and is
$w$-contractive. Since the inequality $w_{k\lambda}(Tx,Ty)\le w_\lambda(x,y)$
with $0<k=L(b_1-a)<1$ holds for all $\lambda>0$, in the iterative procedure
in the proof of Theorem~\ref{t:main} it suffices to choose any $\ov x\in X_w^*$
such that $w_{\ov\lambda}(\ov x,T\ov x)<\infty$ for some $\ov\lambda>0$.
Since $(x_0)'=0$, by virtue of \eq{e:Txx} and \eq{e:abfa}, we find
  $$w_{\lambda_0}(Tx_0,x_0)\le C_0=\frac{b_1-a}{\lambda_0}\,\vfi(L|x_0-y_0|)+
     \frac{\lambda_2}{\lambda_0}\,C_2<\infty$$
(the constants $\lambda_2$ and $C_2$ being evaluated on the interval $[a,b_1]$)
with $\ov\lambda=\lambda_0=L(b_1-a)+1+\lambda_2$, and so, we may set
$\ov x=x_0$. Now, by Theorem~\ref{t:main}, the integral operator $T$ admits
a fixed point: the equality $Tx=x$ on $[a,b_1]$ for some $x\in X_w^*$ is,
by virtue of \eq{e:intop} and \eq{e:aetxt}, equivalent to~\eq{e:ivp}.
\qed\end{pf}

\section{Concluding remarks} \label{s:cr}

\smallbreak
{\bf\ref{s:cr}.1.} It is not our intention in this paper to study the properties
of solutions to \eq{e:ivp} in detail: after Theorem~\ref{t:Cde} on local solutions
of \eq{e:ivp} has been established, the questions of uniqueness, extensions, etc.,
of solutions can be studied following the same pattern as in, e.g., \cite{Filippov}.
Theorems~\ref{t:meniq} and \ref{t:Cde} are valid (with the same proofs) for
mappings $x:[a,b]\to M$ and $f:[a,b]\times M\to M$ satisfying (C.1) and (C.2),
where $(M,|\cdot|)$ is a reflexive Banach space; the details concerning the equality
\eq{e:intab} in this case can be found in \cite{gens}--\cite{Sovae}.

\smallbreak
{\bf\ref{s:cr}.2.} In the theory of the Carath{\'e}odory differential equations
\eq{e:ivp} (cf.\ \cite{Filippov}) the usual assumption on the right hand side is
of the form $|f(t,x)|\le g(t)$ for almost all $t\in[a,b]$ and all $x\in\RB$, where
$g\in\mbox{L}^1[a,b]$, and the resulting solution belongs to $\mbox{AC}[a,b_1]$
for some $a<b_1<b$. However, it is known from \cite[II.8]{KR} that
$\mbox{L}^1[a,b]=\bigcup_{\vfi\in\mathcal{N}}\mbox{L}^\vfi[a,b]$, where
$\mathcal{N}$ is the set of all \mbox{$\vfi$-functions} satisfying the Orlicz condition at
infinity. Also, it follows from \cite[Coro\-lla\-ry~11]{gens} that
$\mbox{AC}[a,b]=\bigcup_{\vfi\in\mathcal{N}}\GV\vfi[a,b]$. Thus,
Theorem~\ref{t:Cde} reflects the {\em regularity\/} property of solutions
of \eq{e:ivp}. Note that, in contrast with functions from $\mbox{AC}[a,b]$,
functions $x$ from $\GV\vfi[a,b]$ have the ``qualified'' modulus of continuity
(\cite[Lemma~3.9(a)]{Sovae}):
$|x(t)-x(s)|\le C_x\cdot\omega_\vfi(|t-s|)$ for all $t,s\in[a,b]$, where
$C_x=d_w^*(x,0)$ and $\omega_\vfi:\RB^+\to\RB^+$ is a subadditive
function given by $\omega_\vfi(u)=u\vfi^{-1}(1/u)$ for $u>0$ and
$\omega_\vfi(+0)=\omega_\vfi(0)=0$.

\smallbreak
{\bf\ref{s:cr}.3.} Theorem~\ref{t:meniq} does not reflect all the flavour of
Theorem~\ref{t:main}, namely, the {\em locality\/} of condition \eq{e:mc}
and the {\em modular convergence\/} of the successive approximations
of the fixed points, and so, an appropriate example is yet to be found;
however, one may try to adjust Example~2.15 from \cite{Kh-Ko-Re}
(note that Proposition~2.14 from \cite{Kh-Ko-Re} is similar to our assertion
\eq{e:corLip} with $k=1$).

\smallbreak
{\em Acknowledgments.} The author is grateful to Marek Balcerzak and Jacek Jachymski
({\L}{\'o}dz Technical University, {\L}{\'o}dz, Poland) for their keen interest and
stimulating discussions on the results of this paper.

\end{document}